\documentstyle[amstex,amssymb]{article}

\renewcommand{\epsilon}{\varepsilon}
\renewcommand{\phi}{\varphi}

\newtheorem{theorem}{Theorem}[section]
\newtheorem{proposition}[theorem]{Proposition}
\newtheorem{lemma}[theorem]{Lemma}
\newtheorem{corollary}[theorem]{Corollary}
\newtheorem{remark}[theorem]{Remark}
\newtheorem{definition}[theorem]{Definition}
\newtheorem{question}[theorem]{Question}
\newtheorem{example}[theorem]{Example}
\newtheorem{conjecture}[theorem]{Conjecture}

\begin{document}
\thispagestyle{empty}

\begin{center}
{\Large\bf On Quasidiagonal $C^*$-algebras}
\end{center}

\vspace{5mm}

\begin{center}{\bf Nathanial P. Brown \footnote{ 
      NSF Postdoctoral Fellow.}}\\
      UC-Berkeley
\end{center}

\begin{abstract}
We give a detailed survey of the theory of quasidiagonal
$C^*$-algebras.  The main structural results are presented and various
functorial questions around quasidiagonality are  discussed.  In
particular we look at what is currently known (and not known) about
tensor products, quotients, extensions, free products, etc.\ of
quasidiagonal $C^*$-algebras.  We also point out how quasidiagonality
is connected to some important open problems.  
\end{abstract}

\parskip2mm

\section{Introduction}

Quasidiagonal $C^*$-algebras have now been studied for more than 20
years.  They are a large class of algebras which arise naturally in
many contexts and include many of the basic examples of finite
$C^*$-algebras.  Notions around quasidiagonality have also played an
important role in BDF/KK-theory and are connected to some of the most
important open questions in $C^*$-algebras. For example, whether every
nuclear $C^*$-algebra satisfies the Universal Coefficient Theorem,
Elliott's Classification Program and whether or not $Ext( C^*_r ({\Bbb
F}_2))$ is a group.

In these notes we give a detailed survey of the basic theory of
quasidiagonal $C^*$-algebras.  At present there is only one survey
article in the literature which deals with this subject (cf.\ [Vo4]).
While there is certainly overlap between this article and [Vo4], the
focus of the present paper is quite different.  We will spend a fair
amount of time giving detailed proofs of a number of basic facts about
quasidiagonal $C^*$-algebras.  Some of these results have appeared in
print, some are well known to the experts but have not (explicitly)
appeared in print and some of them are new.  Moreover, there have been
a number of important advances since the writing of [Vo4].  We will
not give proofs of most of the more difficult recent results.
However, we have tried to at least give precise statements of these
results and have included an extensive bibliography so that the
interested reader may track down the original papers. 

In this paper we are primarily concerned with basic structural
questions.  In particular this means that many interesting topics have
been left out or only briefly touched upon.  For example, we do not
explore the connections between (relative) quasidiagonality and
BDF/KK-theory (found in the work of Salinas, Kirchberg, S. Wassermann,
D\u{a}d\u{a}rlat-Eilers and others) or the generalized inductive limit
approach (introduced by Blackadar and Kirchberg).  But, for the
interested reader, we have included a section containing references to
a number of these topics.  

Throughout the main body of these notes we will only be concerned with
separable $C^*$-algebras and representations on separable Hilbert
spaces.  It turns out that one can usually reduce to this case 
so we don't view this as a major loss of generality.  However, it
causes one problem in that certain nonseparable $C^*$-algebras
naturally arise in the (separable) theory.  Hence we have included an
appendix which deals with the nonseparable case.  

A brief overview of this paper is as follows.  In section 2 we collect
a number of facts that will be needed in the rest of these notes.
This is an attempt to keep the paper self contained, but these
results include some of the deepest and most important tools in
$C^*$-algebra theory and no proofs of well known results are given.

Section 3 contains the definitions and some basic properties of
quasidiagonal operators, quasidiagonal sets of operators and
quasidiagonal (QD) $C^*$-algebras.  We also give some examples
of QD and non-QD $C^*$-algebras.  We end the section with the well
known fact that quasidiagonality implies stable finiteness. 

In section 4 we prove the abstract characterization of QD $C^*$-algebras
which is due to Voiculescu.  

Section 5 deals with the local approximation of QD $C^*$-algebras.  We
show that every such algebra can be locally approximated by a
residually finite dimensional algebra.  We also state a result of
D\u{a}d\u{a}rlat showing that every exact QD $C^*$-algebra can be locally
approximated by finite dimensional $C^*$-algebras.

Section 6 contains the simple fact that every unital QD $C^*$-algebra
has a tracial state.  

Sections 7 - 11 deal with how quasidiagonality behaves under some of
the standard operator algebra constructions.  Section 7 discusses the
easiest of these questions.  Namely what happens when taking
subalgebras, direct products and minimal tensor products of QD
$C^*$-algebras.  Quotients of QD algebras are treated in section 8,
inductive limits in section 9, extensions in section 10 and crossed
products in section 11.  

Section 12 discusses the relationship between quasidiagonality and
nuclearity.  We state a result of Popa which led many experts to
believe that simple QD $C^*$-algebras with 'sufficiently many
projections' are always nuclear.  We then state a result of D\u{a}d\u{a}rlat
which shows that this is not the case.  We also discuss a certain
converse to this question and it's relationship to the classification
program.  Namely the question (due to Blackadar and Kirchberg) of
whether every nuclear stably finite $C^*$-algebra must be QD. 

Section 13 contains miscellaneous results which didn't quite fit
anywhere else.  We state results of Boca and Voiculescu which concern
full free products and homotopy invariance, respectively.  We observe
that all projective algebras and semiprojective MF algebras must be
residually finite dimensional. Finally, we discuss how
quasidiagonality relates to the question of when the classical BDF
$Ext( \cdot )$ semigroups are actually groups and the question of
whether all nuclear $C^*$-algebras satisfy the Universal Coefficient
Theorem.
 
In section 14 we point out where the interested reader can go to learn
more about some of the things that are not covered thoroughly here.

Finally at the end we have an appendix which treats the case of
nonseparable QD $C^*$-algebras.  The main result being that a
$C^*$-algebra is QD if only if all of it's separable $C^*$-subalgebras
are QD.


\section{Preliminaries}

Central to much of what will follow is the theory of {\em completely
positive} maps.  We refer the reader to [Pa] for a comprehensive
treatment of these important maps.  Perhaps the single most important
result about these maps is Stinespring's Dilation Theorem (cf.\ [Pa,
Thm.\ 4.1]).  We will not state the most general version; for our
purposes the following result will suffice.

\begin{theorem} (Stinespring)
Let $A$ be a unital separable $C^*$-algebra and $\phi : A \to B(H)$ be
a unital completely positive map.  Then there exists a separable
Hilbert space $K$, an isometry $V : H \to K$ and a unital
representation $\pi : A \to B(K)$ such that $\phi(a) = V^* \pi(a) V$
for all $a \in A$.
\end{theorem}

Throughout most of [Pa], only the unital case is treated.  The
following result shows that this is not a serious problem.

\begin{proposition}(cf.\ [CE2, Lem.\ 3.9])
Let $\phi : A \to B$ be a contractive completely positive map.  Then
the unique unital extension $\tilde{\phi} : \tilde{A} \to \tilde{B}$ is
also completely positive, where $\tilde{A}, \ \tilde{B}$ are the
$C^*$-algebras obtained by adjoining new units.
\end{proposition} 

We should point out, however, that the easiest proof of this innocent
looking result depends on the nonunital version of Stinespring's
theorem, which is due to E.C. Lance, and does not immediately follow
from the unital version of Stinespring's theorem (cf.\ [La1, Thm.\
4.1]).  (A proof of this result which does not depend on Lance's
result can be given, but it is surprisingly technical.  See [SS, pg.\
14-15].)

Another fundamental result concerning completely positive maps is
Arveson's Extension Theorem.  To state the theorem we need the
following definition.

\begin{definition}
{\em Let $A$ be a unital $C^*$-algebra and $X \subset A$ be a closed linear
subspace. Then $X$ is called an {\em operator system} if $1_A \in X$
and $X = X^*$.}
\end{definition}

\begin{theorem} (Arveson's Extension Theorem) 
If $A$ is a unital $C^*$-algebra, $X \subset A$ is an operator system
and $\phi: X \to C$ is a contractive completely positive map with $C =
B(H)$ or $dim(C) < \infty$ then there exists a completely positive map
$\Phi: A \to C$ which extends $\phi$ (i.e. $\Phi|_X = \phi$).  If $X$
is a $C^*$-subalgebra of $A$ then there always exists a unital
completely positive extension of $\phi$ (whether or not $X$ contains
the unit of $A$). 
\end{theorem}

A proof of the unital statement above can be found in [Pa] while the
nonunital statement is an easy consequence of the unital statement
together with Proposition 2.2. (The nonunital version can also be found
in [La1, Thm.\ 4.2].)

Representations of quasidiagonal $C^*$-algebras will be important in
what follows and hence we will need Voiculescu's Theorem (cf.\ [Vo1]).
In fact, we will need a number of different versions of this result.
It will be convenient to have Hadwin's formulation in terms of rank.

\begin{definition}
{\em If $T \in B(H)$ then let $rank(T) =
dim(\overline{TH})$.  }
\end{definition}

\begin{theorem} 
Let $A$ be a unital $C^*$-algebra and $\pi_i : A \to B(H_i)$ be unital
*-representations for $i = 1,2$.  Then there exists a net of unitaries
$U_{\lambda} : H_1 \to H_2$ such that $\| \pi_2 (a) - U_{\lambda}\pi_1
(a) U_{\lambda}^* \| \to 0$ for all $a \in A$ if and only if
$rank(\pi_1 (a)) = rank(\pi_2 (a))$ for all $a \in A$.  If $A$ is
nonunital then there exists such a net of unitaries if and only if
$rank(\pi_1 (a)) = rank(\pi_2 (a))$ for all $a \in A$ and $dim(H_1) =
dim(H_2)$.
\end{theorem}

When such unitaries exist we say that $\pi_1$ and $\pi_2$ are {\em
approximately unitarily equivalent}.  When both $A$ and the underlying
Hilbert spaces are separable one can even arrange the stronger
condition that $\pi_2 (a) - U_{n}\pi_1 (a) U_{n}^*$ is a compact
operator for each $a \in A, \ n \in {\Bbb N}$ (of course, we can take
a sequence of unitaries when $A$ is separable).  When this is the case
we say that $\pi_1$ and $\pi_2$ are {\em approximately unitarily
equivalent modulo the compacts}.  A proof of this stronger (in the
separable case) result can be found in [Dav, Thm.\ II.5.8] or a proof
of the general result can be found in [Had1].

It turns out that one can usually reduce to the case of separable
$C^*$-algebras and Hilbert spaces.  In this case, the following
version of Voiculescu's theorem will be convenient (cf.\ [Dav, Cor.\
II.5.5]).

\begin{theorem}
Let $H$ be a separable Hilbert space and $C \subset B(H)$ be a unital
separable $C^*$-algebra such that $1_H \in C$.  Let $\iota : C
\hookrightarrow B(H)$ denote the canonical inclusion and let $\rho : C
\to B(K)$ be any unital representation such that $\rho (C \cap {\cal
K}(H)) = 0$.  Then $\iota$ is approximately unitarily equivalent
modulo the compacts to $\iota \oplus \rho$.
\end{theorem}

We will be particularly interested in the case that $C \cap {\cal
K}(H) = 0$.

\begin{definition}
{\em Let $\pi : A \to B(H)$ be a faithful representation of a
$C^*$-algebra $A$.  Then $\pi$ is called {\em essential} if $\pi(A)$
contains no nonzero finite rank operators.  }
\end{definition}

\begin{corollary}
Let $A$ be a separable $C^*$-algebra and $\pi_i : A \to B(H_i)$ be
faithful essential representations with $H_i$ separable for $i = 1,2$.
If $A$ is unital and both $\pi_1, \ \pi_2$ are unital then $\pi_1$ and
$\pi_2$ are approximately unitarily equivalent modulo the compacts.
If $A$ is nonunital then $\pi_1$ and $\pi_2$ are always approximately
unitarily equivalent modulo the compacts.
\end{corollary}

We will need one more form of Voiculescu's Theorem.  We have not been
able to find the following version written explicitly in the
literature.  However, the main idea is essentially due to Salinas (see
the proofs of [Sa1, Thm.\ 2.9] and [DHS, Thm.\ 4.2]). 

If $A$ is a separable, unital $C^*$-algebra and $\phi : A \to B(H)$
(with $H$ separable and infinite dimensional) is a unital completely
positive map then we say that $\phi$ is a {\em representation modulo
the compacts} if $\pi \circ \phi : A \to Q(H)$ is a *-homomorphism,
where $\pi$ is the quotient map onto the Calkin algebra.  If $\pi\circ
\phi$ is injective then we say that $\phi$ is a {\em faithful
representation modulo the compacts}.  In this situation we define
constants $\eta_{\phi} (a)$ by $$\eta_{\phi} (a) = 2\max( \| \phi(a^*
a) - \phi(a^*) \phi(a) \|^{1/2}, \| \phi(a a^*) - \phi(a) \phi(a^*)
\|^{1/2} )$$ for every $a \in A$.

\begin{theorem}
Let $A$ be a separable, unital $C^*$-algebra and $\phi : A \to B(H)$
be a faithful representation modulo the compacts.  If $\sigma : A \to
B(K)$ is any faithful, unital, essential representation then there
exist unitaries $U_n : H \to K$ such that $$\limsup\limits_{n \to
\infty} \| \sigma(a) - U_n \phi(a) U_n^* \| \leq \eta_{\phi} (a)$$ for
every $a \in A$.
\end{theorem}
{\noindent\it Proof.} Note that by Corollary 2.9 it suffices to show
that {\em there exists} a representation $\sigma$ satisfying the
conclusion of the theorem since all such representations are
approximately unitarily equivalent. 

Let $\rho : A \to B(L)$ be the Stinespring dilation of $\phi$; i.e.\
$\rho$ is a unital representation of $A$ and there exists an isometry $V : H
\to L$ such that $\phi(a) = V^* \rho(a) V$, for all $a \in A$.  Let $P
= VV^* \in B(L)$ and $P^{\perp} = 1_L - P$.  We claim that for every
$a \in A$, $$\| P^{\perp}\rho(a) P \| \leq \| \phi(a^* a) - \phi(a^*)
\phi(a) \|^{1/2}.$$ This follows from a simple calculation: 
\begin{multline*}
  \begin{aligned} 
     ( P^{\perp}\rho(a) P)^*( P^{\perp}\rho(a) P)
       & =  P\rho(a^*)P^{\perp} \rho(a) P \\[2mm]
       & =  V V^* \rho(a^*a)V V^* - V V^* \rho(a^*)V V^*\rho(a)V V^*  \\[2mm] 
       & =  V\big(\phi(a^*a) - \phi(a^*)\phi(a)\big)V^*.
  \end{aligned}
\end{multline*}

Now write $L = PL \oplus P^{\perp}L$ and decompose the
representation $\rho$ accordingly.  That is, consider the matrix
decomposition $$\rho (a) =
\left(\begin{array}{cc}
\rho (a)_{11} & \rho (a)_{12} \\
\rho (a)_{21} & \rho (a)_{22}
\end{array}\right), $$  
where $\rho (a)_{21} = P^{\perp}\rho(a) P$ and $\rho (a)_{12} = \rho
(a^*)_{21}^*$.  Hence the norm of the matrix
$$\left(\begin{array}{cc}
0 & \rho (a)_{12} \\
\rho (a)_{21} & 0
\end{array}\right)$$
is bounded above by $1/2 \eta_{\phi}(a)$, since 
$\| P^{\perp}\rho(a) P \| \leq \| \phi(a^* a) - \phi(a^*)
\phi(a) \|^{1/2}$.  

Now comes the trick.  We consider the Hilbert space $P^{\perp}L \oplus
PL$ and the representation $\rho^{\prime} : A \to B(P^{\perp}L \oplus
PL)$ given in matrix form as
$$\rho^{\prime} (a) =
\left(\begin{array}{cc}
\rho (a)_{22} & \rho (a)_{21} \\
\rho (a)_{12} & \rho (a)_{11}
\end{array}\right).$$ 
Now using the obvious identification of the Hilbert spaces $$PL \oplus
\big(\bigoplus_{\Bbb N} P^{\perp}L \oplus PL \big) \ \ {\rm and} \ \
\bigoplus_{\Bbb N} L = \bigoplus_{\Bbb N}( PL \oplus P^{\perp}L)$$ a
standard calculation shows that $$\| \rho(a)_{11} \oplus
\rho^{\prime\infty}(a) - \rho^{\infty}(a)\| \leq \eta_{\phi}(a)$$ for
all $a \in A$, where $\rho^{\prime\infty} = \oplus_{\Bbb N}
\rho^{\prime}$ and $\rho^{\infty} = \oplus_{\Bbb N}
\rho$.  Note also that $\rho(a)_{11} = V \phi(a) V^*$.

Now, let $C$ be the linear space $\phi (A) + {\cal K}(H)$.  Note that
$C$ is actually a separable, unital $C^*$-subalgebra of $B(H)$ with
$\pi (C) = A$ where $\pi: B(H) \to Q(H)$ is the quotient map onto the
Calkin algebra.  By Theorem 2.7 we have that $\iota \oplus
\rho^{\prime\infty} \circ \pi$ is approximately unitarily equivalent
modulo the compacts to $\iota$, where $\iota : C \hookrightarrow B(H)$
is the inclusion.  Let $W_n : H \to H \oplus (\oplus_{\Bbb N}
(P^{\perp}L \oplus PL))$ be unitaries such that $$\| \phi (a) \oplus
\rho^{\prime\infty} (a) - W_n \phi (a) W_n^* \| \to 0$$ for all $a \in
A$.

We now let $$\tilde{V} : H \oplus (\bigoplus_{\Bbb N} (P^{\perp}L
\oplus PL )) \to \bigoplus_{\Bbb N} L$$ be the unitary $V \oplus 1$
(again using the obvious identification of $PL \oplus
\big(\oplus_{\Bbb N} (P^{\perp}L \oplus PL) \big)$ and $\oplus_{\Bbb
N} L$).  Note that $\tilde{V}(\phi(a) \oplus
\rho^{\prime\infty}(a))\tilde{V}^* = V\phi(a)V^* \oplus
\rho^{\prime\infty}(a) = \rho(a)_{11} \oplus
\rho^{\prime\infty}(a)$. We now complete the proof by defining $$K =
\oplus_{\Bbb N} L, \ \ \ \sigma = \rho^{\infty} = \oplus_{\Bbb N}
\rho, \ \ \ U_n = \tilde{V}W_n : H \to \oplus_{\Bbb N} L = K. \ \ \
\Box$$

Finally, we will need a basic result concerning quotient maps of
locally reflexive $C^*$-algebras.  The notion of local reflexivity in
the category of $C^*$-algebras was first introduced by Effros and
Haagerup (cf.\ [EH]). 

\begin{definition} 
{\em A unital $C^*$-algebra $A$ is called {\em locally reflexive} if
each unital completely positive map $\phi : X \to A^{**}$ is the limit
(in the point-weak$^*$ topology) of a net of unital completely
positive maps $\phi_{\lambda} : X \to A$, where $X$ is an arbitrary
finite dimensional operator system and $A^{**}$ denotes the enveloping
von Neumann algebra of $A$.  }
\end{definition}

\begin{definition}
{\em Let $\pi : A \to B$ be a surjective *-homomorphism with $A$
unital.  Then $\pi$ is called {\em locally liftable} if for each
finite dimensional operator system $X \subset B$ there exists a unital
completely positive map $\phi : X \to A$ such that $\pi \circ \phi =
id_X$.}
\end{definition}
  
Of course, if either $A$ or $B$ is {\em nuclear} then the Choi-Effros
Lifting Theorem (cf.\ [CE2, Thm.\ 3.10]) implies that $\pi$ is more than
just locally liftable; there then exists a completely positive
splitting defined on all of $B$.  However, local liftability is all we
will need for our results and it is precisely the class of locally
reflexive $C^*$-algebras which always has this property.  The
following result is a consequence of [EH; 3.2, 5.1, 5.3 and 5.5].

\begin{theorem}
Let $0 \to I \stackrel{\iota}{\to} E \stackrel{\pi}{\to} B \to 0$ be
an exact sequence with $E$ unital.  Then $E$ is locally reflexive if
and only if both $I$ and $B$ are locally reflexive and the morphism
$\pi$ is locally liftable.
\end{theorem}

Local reflexivity plays an important role in the theory of operator
spaces.  We will not need any more results about local reflexivity.
However, we do wish to point out the following implications: $$Nuclear
\Longrightarrow Exact \Longrightarrow Locally \ Reflexive.$$ These
results (together with the definitions of nuclear and exact
$C^*$-algebras) can essentially be found in S. Wassermann's monograph
[Wa3]. ([Wa3] Propositions 5.5 and 5.4 give the first implication
while [Wa3, Remark 9.5.2] states that exactness is equivalent to
property C of Archbold and Batty.  However, property C implies
property C$^{\prime\prime}$, as defined in [EH, pg.\ 120], which in
turn is equivalent to local reflexivity by [EH, Thm.\ 5.1].)  Since
the pioneering work of E. Kirchberg, exactness has played a central
role in $C^*$-algebras.  However, since we will only need the local
liftability statement of Theorem 2.13, we will also consider the class
of locally reflexive $C^*$-algebras.


\section{Definitions, Basic Results and Examples}

Recall that throughout the body of these notes all Hilbert spaces and
$C^*$-algebras are assumed to be separable.

We begin this section by recalling the notions of {\em block diagonal}
and {\em quasidiagonal} operators on a Hilbert space.  In Proposition
3.4 we show that the notion of a quasidiagonal operator can be
expressed in terms of a local finite dimensional approximation
property. This local version then extends to a suitable definition of
a quasidiagonal (QD) $C^*$-algebra (Definition 3.8).  In Theorem 3.11 we
prove a fundamental result about representations of QD $C^*$-algebras.
At the end of this section we give some examples of QD (and non-QD)
$C^*$-algebras and observe that QD $C^*$-algebras are always {\em
stably finite} (cf.\ Proposition 3.19).

\begin{definition}
{\em A bounded linear operator $D$ on a  Hilbert space $H$
is called {\em block diagonal} if there exists an increasing sequence
of finite rank projections, $P_1 \leq P_2 \leq P_3 \cdots$, such that
$\| [D, P_n] \| = \| DP_n - P_n D \| = 0$ for all $n \in {\Bbb N}$ and
$P_n \to 1_H$ (in the strong operator topology) as $n \to \infty$.}
\end{definition}

Note that if $\| [D, P_n] \| = 0$ then $\| [D, (P_n - P_{n - 1})] \| =
0$ as well.  Thus the matrix for $D$ with respect to the decomposition
$H = P_1 H \oplus (P_2 - P_1)H \oplus (P_3 - P_2)H \oplus \cdots$ is
block diagonal.

The notion of a quasidiagonal operator is due to Halmos and is a
natural generalization of a block diagonal operator.

\begin{definition}
{\em A bounded linear operator $T$ on a  Hilbert space $H$
is called {\em quasidiagonal} if there exists an increasing sequence
of finite rank projections, $P_1 \leq P_2 \leq P_3 \cdots$, such that
$\| [T, P_n] \| = \| TP_n - P_n T \| \to 0$ and $P_n \to 1_H$ (in the
strong operator topology) as $n \to \infty$.}  
\end{definition}

Halmos observed the following relationship between quasidiagonal and
block diagonal operators.  

\begin{proposition}
If $T \in B(H)$ then $T$ is quasidiagonal if
and only if there exist a block diagonal operator $D \in B(H)$ and a
compact operator $K \in {\cal K}(H)$ such that $T = D + K$.
\end{proposition}

We will not give the proof of this proposition here as it is a special
case of Theorem 5.2.  Note, however, that one direction is easy.
Namely, if $T = D + K$ as above then $T$ must be quasidiagonal since
any increasing sequence of finite rank projections converging to $1_H$
(s.o.t.) will form an approximate identity for ${\cal K}(H)$ and hence
will asymptotically commute with every compact operator.

It is an important fact that the seemingly global notion of
quasidiagonality can be expressed in a local way.

\begin{proposition}
Let $T \in B(H)$.  Then $T$ is quasidiagonal if and only if for each
finite set $\chi \subset H$ and $\varepsilon > 0$ there exists a
finite rank projection $P \in B(H)$ such that $\| [T, P] \| \leq
\varepsilon$ and $\| P(x) - x \| \leq \varepsilon$ for all $x \in
\chi$.
\end{proposition} 
{\it\noindent Proof.}  We may assume that $\| T \| \leq 1$. It is
clear that the definition of a quasidiagonal operator implies the
condition stated above.  To prove the converse, it suffices to show
that for each finite set $\chi \subset H$ and $\varepsilon > 0$ there
exists a finite rank projection $P$ such that $\| [P, T] \| <
\varepsilon$ and $P(x) = x$ for all $x \in \chi$.  Having established
this it is not hard to construct finite rank projections $P_1 \leq P_2
\leq P_3 \cdots$, such that $\| [T, P_n] \| \to 0$ and $P_n \to 1_H$
in the strong operator topology.

So let $\chi \subset H$ be a finite set, $\varepsilon > 0$ and let $R$
be the orthogonal projection onto $K = span\{ \chi \}$.  By
compactness of the unit ball of $K$ there is a finite set
$\tilde{\chi} \subset K$ which is $\varepsilon$-dense in the unit ball
of $K$.  Now let $Q$ be a finite rank projection such that $\| [Q, T]
\| < \varepsilon$ and $\| Q(x) - x \| < \varepsilon$ for all $x \in
\tilde{\chi}$.  Then for all $y \in K$ we have $\| Q(y) - y \| <
3\varepsilon \| y \|$ and hence $\| (1 - R)QR \| < 3 \varepsilon$. 

Now consider the positive contraction $X = RQR + (1 - R)Q(1 - R)$.
Observe that $X$ is actually very close to $Q$:
\begin{multline*}
  \begin{aligned}
    \|Q - X \| 
    & = \|RQ(1 - R) + (1 - R)QR \| \\[2mm]
    & = max \{ \ \| RQ(1 - R) \|, \ \| (1 - R)QR  \| \ \}\\[2mm] 
    & = \| (1 - R)QR \| \\[2mm] 
    & < 3\varepsilon.
   \end{aligned}
\end{multline*}
Hence $X$ is almost a projection (i.e.\ it's spectrum is contained in
$[0, 3\varepsilon) \cup (1 - 3\varepsilon, 1]$).  Let $P$ be the
projection obtained from functional calculus on $X$.  Then $\| P - Q
\| \leq 6\varepsilon$ and hence $\| [P, T] \| \leq 13\varepsilon$.
Finally we claim that $P(x) = x$ for all $x \in \chi$.  To see this,
first note that $X$ commutes with $R$ and hence so does $P$.  This
implies that $PR = RPR$ is a projection with support contained in $K$.
However, for each $y \in K$ we also have $\| PR(y) - y \| = \| R(P(y)
- y) \| \leq \| P(y) - Q(y) \| + \| Q(y) - y \| \leq 9\varepsilon
\|y\|.$ Hence the support of $PR$ is all of $K$; i.e. $PR = R$. $\Box$

With this local characterization in hand we now define the following
generalization of a quasidiagonal operator.  

\begin{definition}
{\em A subset $\Omega \subset B(H)$ is a called a {\em quasidiagonal
set of operators} if for each finite set $\omega \subset \Omega$,
finite set $\chi \subset H$ and $\varepsilon > 0$ there exists a
finite rank projection $P \in B(H)$ such that $\| [T, P] \| \leq
\varepsilon$ and $\| P(x) - x \| \leq \varepsilon$ for all $T \in
\omega$ and $x \in \chi$.}
\end{definition}

It is easy to see that a set $\Omega \subset B(H)$ is a quasidiagonal
set of operators if and only if the $C^*$-algebra generated by
$\Omega$, $C^*(\Omega) \subset B(H)$, is a quasidiagonal set of
operators.

The proof of the next proposition is a straightforward adaptation of
the proof of Proposition 3.4.

\begin{proposition}
If $A \subset B(H)$ is separable then $A$ is a quasidiagonal set of
operators if and only if there exists an increasing sequence of finite
rank projections, $P_1 \leq P_2 \leq P_3 \cdots$, such that for all $a
\in A$, $\| [a, P_n] \| \to 0$ and $P_n \to 1_H$ (s.o.t.) as $n
\to \infty$.
\end{proposition}

\begin{remark}
The previous proposition is often used when defining quasidiagonal
$C^*$-algebras.  However, L. Brown has pointed out to us that the
previous proposition is not true if $A$ is not separable (even if $H$
is separable).  Definition 3.5 allows one to use Zorn's lemma to
construct {\em maximal} quasidiagonal subsets of $B(H)$ and we claim
that they provide counterexamples.  The proof goes by contradiction.
So assume that $\Omega \subset B(H)$ is a maximal quasidiagonal set of
operators and there exist finite rank projections such that $\| [x,
P_n] \| \to 0$ for all $x \in \Omega$.  Construct a block diagonal
operator $T$ such that $[T, P_{2n}] = 0$ and $\| [T, P_{2n + 1}] \| =
1$ for all $n \in {\Bbb N}$.  Since $[T, P_{2n}] = 0$ for all $n$ we
see that $\Omega \cup \{ T \}$ is a quasidiagonal set of operators and
hence (by maximality) $T \in \Omega$.  This gives the contradiction
since $\| [T, P_{2n + 1}] \| = 1$ for all $n \in {\Bbb N}$.  Hence it
is important to take Definition 3.5 for general quasidiagonal
questions.
\end{remark}

We finally come to the definition of a quasidiagonal $C^*$-algebra.

\begin{definition}
{\em Let $A$ be a $C^*$-algebra.  Then $A$ is called {\em
quasidiagonal} (QD) if there exists a faithful representation $\pi : A
\to B(H)$ such that $\pi(A)$ is a quasidiagonal set of operators.}
\end{definition}

Some remarks regarding this definition are in order.  First we should
point out that some authors (e.g.\ [Had2]) refer to $C^*$-algebras
satisfying Definition 3.8 as `weakly quasidiagonal' $C^*$-algebras.
There is good reason for this terminology as it emphasizes the
distinction between abstract and concrete $C^*$-algebras.  It is
important to make this distinction since every $C^*$-algebra has a
representation $\pi$ such that $\pi (A)$ is a quasidiagonal set of
operators (namely the zero representation).  On the other hand, it is
possible to give examples of $C^*$-algebras $A$ and faithful
representations $\pi : A \to B(H)$ such that $A$ is QD but $\pi(A)$ is
{\em not} a quasidiagonal set of operators.  Perhaps the most extreme
case of this is an example of L. Brown.  In [BrL2] it was shown that
there exists an operator $T$ on a separable Hilbert space such that
$T\oplus T$ is quasidiagonal while $T$ is not!  Hence $C^*(T)$ is a QD
$C^*$-algebra but is {\em not} a quasidiagonal set of operators in
it's natural representation.  Thus it is indeed very important to
distinguish between abstract QD $C^*$-algebras and concrete
quasidiagonal sets of operators. (The reader is cautioned that this is
not always done carefully in the literature.)  Other authors prefer to
say that a {\em representation} is quasidiagonal if it's image is a
quasidiagonal set of operators.  Definition 3.8 then becomes
equivalent to the statement that $A$ admits a faithful quasidiagonal
representation.

Definitions 3.5 and 3.8 are the correct definitions in the nonseparable
case as well (see the Appendix).  We will see that certain nonseparable
$C^*$-algebras (namely $\Pi M_n ({\Bbb C})$) naturally arise in the
separable theory and hence it will be logically necessary to treat
this case also.

Finally note that Definition 3.8 does not require the representation
to be nondegenerate.  Of course, this can always be arranged.  Note,
however, that this is actually a deep fact as the proof of Lemma 3.10
below depends on Voiculescu's Theorem (at least in the nonunital
case).

\begin{definition}
{\em If $\pi : A \to B(H)$ is a representation and $L \subset H$ is a
$\pi(A)$-invariant subspace then $\pi_L : A \to B(L)$ denotes the
restriction representation (i.e.\ $\pi_L (a) = P_L \pi(a)|_L$, where
$P_L$ is the orthogonal projection from $H \to L$). }
\end{definition}

\begin{lemma}
Let $\pi: A \to B(H)$ be a faithful representation and $L \subset H$
be the nondegeneracy subspace of $\pi(A)$.  Then $\pi(A)$ is a
quasidiagonal set of operators if and only if $\pi_L (A)$ is a
quasidiagonal set of operators.
\end{lemma}
{\it\noindent Proof.}  Assume first that $\pi_L (A)$ is a
quasidiagonal set of operators.  Then write $H = L \oplus \tilde{L}$.
Since $\pi(a) = \pi_L (a) \oplus 0$, any finite rank projection $P \in
B(L)$ can be extended to a finite rank projection $P \oplus \tilde{P}$
such that $\| [\pi(a), P \oplus \tilde{P}] \| = \| [\pi_L(a), P] \|$.
From this one deduces that $\pi(A)$ must also be a quasidiagonal set
of operators.

Now assume that $\pi(A)$ is a quasidiagonal set of operators.  If $A$
is unital then $\pi(1_A) = P_L$.  If $R \in B(H)$ is any finite rank
projection that almost commutes with $\pi(1_A) = P_L$ then $P_L R P_L
\in B(L)$ is very close to a projection which does commute with $P_L$.
We leave the details to the reader, but some standard functional
calculus then implies that $\pi_L (A)$ is also a quasidiagonal set of
operators.

In the case that $\pi(A)$ is a quasidiagonal set of operators and $A$
is nonunital, we have to call on Voiculescu's Theorem (version 2.6).
Since it is clear that $rank(\pi(a)) = rank(\pi_L (a))$ for all $a \in
A$ we have that $\pi$ and $\pi_L$ are approximately unitarily
equivalent.  However, it is an easy exercise to verify that if $\rho$
and $\tilde{\rho}$ are two approximately unitarily equivalent
representations then $\rho(A)$ is a quasidiagonal set of operators if
and only if $\tilde{\rho} (A)$ is a quasidiagonal set of operators.
$\Box$

We now give the fundamental theorem about representations of QD
$C^*$-algebras.

\begin{theorem}(cf.\ [Vo4, 1.7])
Let $\pi : A \to B(H)$ be a faithful {\em essential} (cf.\ Definition
2.8) representation.  Then $A$ is QD if and only if $\pi (A)$ is a
quasidiagonal set of operators.
\end{theorem}
{\it\noindent Proof.} If $\pi(A)$ is a quasidiagonal set of operators
then, of course, $A$ is QD.  Conversely, if $A$ is QD then there
exists a faithful representation $\rho : A \to B(K)$ such that
$\rho(A)$ is a quasidiagonal set of operators.  In light of Lemma
3.10, we may assume that both $\pi$ and $\rho$ are
nondegenerate. Defining $\rho_{\infty} = \oplus_{\Bbb N} \rho : A \to
B(\oplus_{\Bbb N} K)$ it is easy to see that $\rho_{\infty} (A)$ is
also a quasidiagonal set of operators.  But, since $\rho_{\infty}$ is
an essential representation, Voiculescu's Theorem (version 2.9)
implies that $\pi$ and $\rho_{\infty}$ are approximately unitarily
equivalent.  Hence $\pi(A)$ is also a quasidiagonal set of operators.
$\Box$

We now give some examples of QD $C^*$-algebras and non-QD $C^*$-algebras.

\begin{example}
Every commutative $C^*$-algebra is QD.  Indeed, if $A = C_0 (X)$ and
for each $x \in X$ we let $ev_x : A \to {\Bbb C}$ be evaluation at $x$
then $\pi = \oplus_{x \in F} ev_x$, where $F \subset X$ is a countable
dense set, is a faithful representation and it is easy to see that
$\pi (A)$ is a quasidiagonal (in fact, diagonal) set of operators.
\end{example}

\begin{example}
Approximately finite dimensional (AF) algebras are QD.  Let $A =
\overline{\cup_n A_n}$ be AF with each $A_n \subset A_{n + 1}$ finite
dimensional.  Let $\pi : A \to B(H)$ be a faithful nondegenerate
representation and write $H = \overline{\cup_n H_n}$ where each $H_n
\subset H_{n + 1}$ is a finite dimensional subspace.  Then define $P_n
\in B(H)$ to be the (finite rank) projection onto the subspace $\pi
(A_n)H_n$.  Then we evidently have that $\| [\pi (a), P_n] \| \to 0$
for all $a \in A$ and $P_n \to 1_H$ in the strong topology.
\end{example}

\begin{example}
Irrational rotation algebras are QD.  That is, if $A_{\theta}$ is the
universal $C^*$-algebra generated by two unitaries $U, V$ subject to
the relation $UV = (\exp (2\pi \theta i)) VU$ for some irrational number
$\theta \in [0, 1]$ then $A_{\theta}$ is QD.  This was first proved by
Pimsner and Voiculescu when they showed how to embed $A_{\theta}$ into
an AF algebra (cf.\ [PV2]).  This was later generalized by Pimsner in
[Pi] (see also section 11 of these notes).  
\end{example}

\begin{example}
Perhaps the most important class of QD $C^*$-algebras are the
so-called residually finite dimensional (RFD) $C^*$-algebras.  A
$C^*$-algebra $R$ is called RFD if for each $x \in R$ there exists a
*-homomorphism $\pi : R \to B$ such that $dim(B) < \infty$ and $\pi(x)
\neq 0$.  That such algebras have a faithful representation whose
image is a quasidiagonal (in fact, block diagonal) set of operators is
proved similar to the case of abelian algebras.  Often times general 
questions about QD $C^*$-algebras can be reduced to the case of RFD algebras.
\end{example}

\begin{example}
Both the cone ($CA = C_0 ((0,1]) \otimes A$) and suspension ($SA = C_0
((0,1)) \otimes A$) over any $C^*$-algebra $A$ are QD.  Since $SA
\subset CA$ and $CA$ is homotopic to $\{0\}$, this can be deduced from
the homotopy invariance of quasidiagonality (cf.\ [Vo3] or Theorem 13.1
of these notes). From this we see that every $C^*$-algebra is a
quotient of a QD $C^*$-algebra (since $A \cong CA/SA$).
\end{example}

\begin{example}
A $C^*$-algebra which contains a proper (i.e.\ non-unitary) isometry
is not QD.  Since it is clear that a subalgebra of a QD $C^*$-algebra
is again QD, it suffices to show that the Toeplitz algebra is not
QD. (Recall that Coburn's Theorem states that the $C^*$-algebras
generated by any two proper isometries are isomorphic.) We let
$C^*(S)$ denote the Toeplitz algebra, where $S$ is a proper isometry,
and let $\pi : C^*(S) \to B(H)$ be any faithful unital essential
representation.  Then $\pi(S)$ is a semi-Fredholm operator with index
$-\infty$.  On the other hand, it follows from Proposition 3.3 that any
semi-Fredholm quasidiagonal operator on $H$ must have index zero
(since any semi-Fredholm block diagonal operator must have index zero)
and hence $\pi(S)$ is not a quasidiagonal operator.  Hence, by Theorem
3.11, $C^*(S)$ is not QD.  (See [Hal1] for generalizations of this
result.)
\end{example}

The previous example implies a more general result.

\begin{definition}
{\em Let $A$ be a unital $C^*$-algebra.  Then $A$ is said to be {\em
stably finite} if $A \otimes M_n ({\Bbb C})$ contains no proper
isometries for all $n \in {\Bbb N}$.  If $A$ is nonunital, then $A$ is
called {\em stably finite} if the unitization $\tilde{A}$ is stably
finite.  }
\end{definition}

\begin{proposition}
QD $C^*$-algebras are stably finite.
\end{proposition}
{\it\noindent Proof.}  It is easy to see that if $A$ is nonunital and
QD then the unitization $\tilde{A}$ is also QD.  Furthermore, it is a
good exercise to verify that if $A$ is QD then for all $n \in {\Bbb
N}$, $M_n({\Bbb C}) \otimes A$ is also QD.  From these observations
and Example 3.17 we see that if $A$ is QD then $M_n({\Bbb C}) \otimes A$
(or $M_n({\Bbb C}) \otimes \tilde{A}$ in the non-unital case) has no
proper isometries for all $n \in {\Bbb N}$.  Hence $A$ is stable
finite.  $\Box$

The converse is not true.  S.\ Wassermann has given examples of non-QD
MF algebras (cf.\ Definition 9.1 and Example 8.6 of these notes).  But
every MF algebras is stably finite (cf.\ [BK1, Prop.\ 3.3.8]).  Hence,
in general, QD is not equivalent to stably finite.  However, Blackadar
and Kirchberg have asked whether or not they are equivalent within the
category of nuclear $C^*$-algebras (see Question 12.5).


\section{Voiculescu's Abstract Characterization}

In this section we prove an abstract (i.e.\ representation free)
characterization of QD $C^*$-algebras (cf.\ [Vo3, Thm.\ 1]).  This
fundamental result will be crucial in sections 8 - 10.

Consider the following property of an arbitrary $C^*$-algebra $A$.

($\ast$) For each finite set ${\cal F} \subset A$ and $\varepsilon >
0$ there exists a contractive completely positive map $\phi : A \to B$
such that $i) \ dim(B) < \infty$, $ii) \ \|\phi (x) \| \geq \| x \| -
\varepsilon$ for all $x \in {\cal F}$ and $iii) \ \| \phi (xy) -
\phi(x)\phi(y) \| \leq \varepsilon$ for all $x,y \in {\cal F}$. 

For a unital algebra $A$ we have a related property.

($\ast\ast$) For each finite set ${\cal F} \subset A$ and $\varepsilon > 
0$ there exists a unital completely positive map $\phi : A \to B$ such
that $i) \ B \cong M_n({\Bbb C})$ for some $n \in {\Bbb N}$, $ii) \
\|\phi (x) \| \geq \| x \| - \varepsilon$ for all $x \in {\cal F}$ and
$iii) \ \| \phi (xy) - \phi(x)\phi(y) \| \leq \varepsilon$ for all $x,y
\in {\cal F}$.

We will refer to such maps as {\em $\epsilon$-isometric} and {\em
$\epsilon$-multiplicative} on ${\cal F}$.

\begin{lemma}
If $A$ is a unital $C^*$-algebra then $A$ satisfies ($\ast$) if and
only if $A$ satisfies ($\ast\ast$).
\end{lemma}
{\it\noindent Proof.} ($\Leftarrow$) This is obvious.

($\Rightarrow$) We only sketch the main idea.  First, we identify $B$
with a unital subalgebra of $M_m({\Bbb C}) = B({\Bbb C}^m)$ for some
$m \in {\Bbb N}$. Let $1_A$ denote the unit of $A$ and let $P \in
M_m({\Bbb C})$ be the projection onto the range of $\phi(1_A)$.  Then
one shows that $\phi(a) = P\phi(a) = \phi(a) P$ for all $a \in A$.
Moreover, if $\phi$ is very multiplicative on $1_A$ then $\phi(1_A)$
is close to $P$.

Now let $\psi : A \to PM_m({\Bbb C})P \cong M_n({\Bbb C})$ (for some
$n \leq m$) be given by $\psi(a) = P\phi(a)P$ and clearly $\psi$ has
the same multiplicativity and isometric properties (up to
$\varepsilon$) that $\phi$ does.  Moreover, since $\phi(1_A)$ is close
to $P$, $\psi(1_A)$ is invertible in $PM_m({\Bbb C})P \cong M_n({\Bbb
C})$.  Thus we replace $\psi$ with the map $a \mapsto
(\psi(1_A))^{-1/2}\psi(a)(\psi(1_A))^{-1/2}$ to get a unital complete
positive map into a matrix algebra.  The multiplicativity and
isometric properties of this new map are not quite as good as those of
$\phi$, but they are good enough.  $\Box$

We are now ready for Voiculescu's abstract characterization of QD
$C^*$-algebras.  Our proof is based on the proof of [DHS, Thm.\
4.2] and, hopefully, is easier to follow than the original.
However, the main ideas are the same.  We have simply isolated the
hard part in Theorem 2.10.

\begin{theorem}[Voiculescu]
Let $A$ be a $C^*$-algebra.  Then $A$ is QD if and only if $A$
satisfies {\em ($\ast$)}.
\end{theorem}
{\noindent\it Proof.}  From Proposition 2.2 it is easy to see that $A$
satisfies ($\ast$) if and only if $\tilde{A}$ satisfies ($\ast$).
Similarly is it clear that $A$ is QD if and only if $\tilde{A}$ is QD
and hence we may assume that $A$ is unital.

($\Rightarrow$) Let $\pi : A \to B(H)$ be a unital faithful essential
representation on a separable Hilbert space.  We can
then find an increasing sequence of finite rank projections, $P_1 \leq
P_2 \leq P_3 \cdots$, such that for all $a \in A$, $\| [\pi (a), P_n]
\| \to 0$ and $P_n \to 1_H$ in the strong topology.  Then for all $n$,
$P_n B(H) P_n$ is isomorphic to a matrix algebra and the unital
completely positive maps $\phi_n : A \to P_n B(H) P_n$, $a \mapsto P_n
\pi(a) P_n$ are easily seen to be asymptotically multiplicative and
isometric.

($\Leftarrow$) By Lemma 4.1 we can find a sequence of unital
completely positive maps $\phi_i : A \to M_{n(i)} ({\Bbb C})$ which
are asymptotically multiplicative and asymptotically isometric.  Let
$$H_m = \bigoplus_{i = m}^{\infty} {\Bbb C}^{n(i)}, \ \ \ \Phi_m =
\bigoplus_{i = m}^{\infty} \phi_i : A \to B(H_m).$$ Evidently each
$\Phi_m$ is a faithful representation modulo the compacts (as in
Theorem 2.10).  Let $\sigma : A \to B(K)$ be any faithful, unital,
essential representation and by Theorem 2.10 we can find unitaries
$U_m : H_m \to K$ such that $\| \sigma(a) - U_m \Phi_m(a) U_m^* \| \to
0$ as $m \to \infty$ for all $a \in A$.  Since it is clear that
$\Phi_m (A)$ is a quasidiagonal (in fact, block diagonal) set of
operators for every $m$ it is easy to see that $\sigma(A)$ is also a
quasidiagonal set of operators and hence $A$ is QD.  $\Box$

In addition to being a very useful tool in establishing the
quasidiagonality of a given $C^*$-algebra this result also shows that
QD $C^*$-algebras are a very natural {\em abstract} class of algebras.
Indeed, this result shows that QD $C^*$-algebras are precisely those
which have `good matrix models' in the sense that all of the relevant
structure (order, adjoints, multiplication, norms) can approximately
be seen in a matrix.

We wish to note a minor generalization which will be useful later on.

\begin{definition}
{\em If $A$ is a unital $C^*$-algebra and ${\cal F} \subset A$ is a finite
set then we will let $X_{{\cal F}{\cal F}}$ denote the smallest
operator system (cf.\ Definition 2.3) containing ${\cal F}$ and $\{ ab
: a,b \in {\cal F} \}$.}
\end{definition}

\begin{definition}
{\em If $ F, B \subset B(H)$ are sets of operators then we say $F$ is
{\em $\epsilon$-contained in $B$} if for each $x \in F$ there exists
$y \in B$ such that $\| x - y \| < \epsilon$.  When this is the case
we write $F \subset^{\epsilon} B$.  }
\end{definition}

\begin{corollary}
Assume that $A$ is unital and for every finite subset ${\cal F}
\subset A$ and $\varepsilon > 0$ there exists a contractive completely
positive map $\phi : X_{{\cal F}{\cal F}} \to B(H)$ such that $\phi$
is $\varepsilon$-isometric and $\varepsilon$-multiplicative on ${\cal
F}$ and $\phi({\cal F})$ is $\epsilon$-contained in a QD $C^*$-algebra
$B \subset B(H)$. Then $A$ is QD.
\end{corollary}
{\it\noindent Proof.}  That $A$ satisfies $(\ast)$ follows from
Arveson's Extension Theorem applied to $\phi$ and to the almost
isometric and multiplicative maps from $B$ to finite dimensional
$C^*$-algebras. $\Box$

\begin{remark}
Note that the hypotheses of the previous corollary can be relaxed
further.  Indeed, one only needs such $\varepsilon$-isometric and
$\varepsilon$-multiplicative maps on a sequence of finite sets which are 
suitably  dense in $A$ (e.g.\ generate a dense *-subalgebra of $A$).
\end{remark}


\section{Local Approximation}

We observe that every QD $C^*$-algebra can be locally approximated by
a residually finite dimensional (RFD) $C^*$-algebra (cf.\ Example 3.15).
The proof is a simple adaptation of Halmos' original proof that every
quasidiagonal operator can be written as a block diagonal operator
plus a compact. We also recall a result of M.\ D\u{a}d\u{a}rlat which gives a
much stronger approximation in the case of exact QD $C^*$-algebras.

\begin{definition}(cf.\ Definition 3.1)
{\em Let $B \subset B(H)$ be a $C^*$-algebra.  Then $B$ is called a
{\em block diagonal} algebra if there exists an increasing sequence of
finite rank projections, $P_1 \leq P_2 \leq P_3 \cdots$, such that $\|
[b, P_n] \| = 0$ for all $b \in B$, $n \in {\Bbb N}$ and $P_n \to 1_H$
(s.o.t.).  }
\end{definition}

It is relatively easy to see that a $C^*$-algebra $R$ is RFD if and
only if there exists a faithful representation $\pi : R \to B(H)$ such
that $\pi (R)$ is a block diagonal algebra.  The next result, which is
well known to the experts, shows that every QD $C^*$-algebra can be locally
approximated by an RFD algebra.

\begin{theorem}
Let $A \subset B(H)$ be a $C^*$-algebra.  Then $A$ is a quasidiagonal
set of operators if and only if for every finite set ${\cal F} \subset
A$ and $\varepsilon > 0$ there exists a block diagonal algebra $B
\subset B(H)$ such that ${\cal F} \subset^{\epsilon} B$ (cf.\
Definition 4.4) and $A + {\cal K}(H) = B + {\cal K}(H)$.
\end{theorem}
{\it\noindent Proof.} Clearly we only have to prove the necessity
since $B + {\cal K}(H)$ is a quasidiagonal set of operators.  Our
proof follows closely the proof of [Ar, Thm.\ 2] where a similar
result is obtained for general quasicentral approximate units.

So let ${\cal F} \subset A$ and $\varepsilon > 0$ be given.  We may
assume that ${\cal F}$ is contained in the unit ball of $A$.  Let
${\cal F}_1 \subset {\cal F}_2 \subset {\cal F}_3 \ldots$ be a
sequence of finite sets such that ${\cal F} \subset {\cal F}_1$ and
whose union is dense in the unit ball of $A$.  Since $A$ is a
quasidiagonal set of operators we can use Proposition 3.6 to find
finite rank projections $P_1 \leq P_2 \ldots $ converging to $1_H$
(strongly) and such that $\| [P_n , a] \| \to 0$ for all $a \in A$.
By passing to a subsequence we may assume that $\| [P_n , a] \| <
\epsilon/(2^n)$ for all $a \in {\cal F}_n$.  Now, let $E_n = P_{n} -
P_{n - 1}$ for $n = 1, \ 2, \ldots$ where $P_0 = 0$.  Note that
$\sum\limits_{n = 1}^{\infty} E_n = 1_H$.

Then one defines completely positive maps $\delta_k : A \to B(H)$ via
the formula $$\delta_k (a) = \sum\limits_{n = 1}^{k} E_n a E_n.$$ We
leave it to the reader to verify that the $\delta_k$'s converge in the
point strong operator topology (i.e.\ $\delta_k (a)$ is strongly
convergent for each $a \in A$) and hence $$\delta (a) = \sum\limits_{n
= 1}^{\infty} E_n a E_n.$$ is a well defined completely positive map.
Now let $B = C^*(\delta(A))$ and clearly $B$ is a block diagonal set
of operators.  Moreover, for each $a \in A$ we have
\begin{multline*}    
     \begin{aligned}
         a - \delta(a)
          & =  \sum\limits_{n = 1}^{\infty} aE_n - 
               \sum\limits_{n = 1}^{\infty} E_n a E_n  \\[2mm]      
          & =  \sum\limits_{n = 1}^{\infty} (aE_n - E_n a E_n) \\[2mm]
          & =  \sum\limits_{n = 1}^{\infty} (aE_n - E_n a)E_n 
  \end{aligned}
\end{multline*} 
where convergence of these sums is again taken in the strong operator
topology.  However, for each $a \in \cup {\cal F}_n$ the last
summation above is actually convergent in the norm topology and is
compact since the $E_n$'s are finite rank.  Note that by construction
we have $\| a - \delta (a) \| \leq \sum \epsilon/(2^n) = \epsilon$ for
all $a \in {\cal F}_1$.  Now since $\delta$ is norm continuous (being
completely positive) we then conclude that $a - \delta(a)$ is a
compact operator for all $a \in A$.  It follows that $A +
{\cal K}(H) = B + {\cal K}(H)$. $\Box$

Theorem 5.2 fails when $A$ is not separable (cf.\ Remark 3.7).

\begin{corollary}(cf.\ [GM])
Every (separable) $C^*$-algebra $A$ is a quotient of an RFD algebra.
If $A$ is nuclear (resp.\ exact) then the RFD algebra can be chosen
nuclear (resp.\ exact).
\end{corollary}
{\noindent\it Proof.} Let $\pi : CA \to B(H)$ be a faithful essential
representation of the cone over $A$ (cf.\ Example 3.16).  If $A$ is
nuclear (resp.\ exact) then so is $CA$ and hence so is $\pi(CA) +
{\cal K}(H)$ (cf.\ [CE1, Cor.\ 3.3], [Kir2, Prop.\ 7.1]).  Let $R
\subset B(H)$ be an RFD algebra such that $\pi(CA) + {\cal K}(H) = R +
{\cal K}(H)$.  Passing to the Calkin algebra we see that $CA$, and
hence $A$, is a quotient of $R$.  Since exactness passes to
subalgebras ([Kir2, Prop.\ 7.1]), it is clear that $R$ is exact
whenever $A$ is exact.  When $A$ is nuclear we deduce that $R$ is also
nuclear from [CE1, Cor.\ 3.3] and the exact sequence $$0 \to R\cap
{\cal K}(H) \to R \to CA \to 0,$$ since ${\cal K}(H)$ is type I and
hence all of it's subalgebras are nuclear (cf.\ [Bl1]).  $\Box$

The next result of D\u{a}d\u{a}rlat is a vast improvement under the
additional assumption of exactness.  We will not prove this here; see
[D\u{a}d3, Thm.\ 6].  However we remark that the proof depends in an
essential way on Theorem 5.2 as it allows one to reduce to the case of
RFD algebras.

\begin{theorem}[D\u{a}d\u{a}rlat] 
Let $A \subset B(H)$ be such that $A \cap {\cal K}(H) = 0$.  Then $A$
is exact and QD if and only if for every finite set ${\cal F} \subset
A$ and $\varepsilon > 0$ there exists a finite dimensional subalgebra
$B \subset B(H)$ such that ${\cal F} \subset^{\epsilon} B$.
\end{theorem}

Note the similarity with the definition of an AF algebra.  The
difference, of course, is that we have had to go outside the algebra
to get the finite dimensional approximation.  We regard this as very
strong evidence in favor of an affirmative answer to the following
conjecture. (See also [BK1, Question 7.3.3])

\begin{conjecture}
Every (separable) exact QD $C^*$-algebra is isomorphic to a subalgebra
of an AF algebra.
\end{conjecture}


\section{Traces}

\begin{proposition}(cf.\ [Vo4, 2.4])
If $A$ is a unital QD $C^*$-algebra then $A$ has a tracial state.
\end{proposition}
{\it\noindent Proof.}  By Theorem 4.2 and Lemma 4.1 we can find a
sequence of unital completely positive maps $\phi_i : A \to M_{n(i)}
({\Bbb C})$ such that $\| a \| = \lim_{i} \| \phi_{i} (a) \|$ and $\|
\phi_{i} (ab) - \phi_{i} (a) \phi_{i} (b) \| \to 0$ for all $a,b \in
A$.  Let $\tau_{n(i)}$ denote the tracial state on $M_{n(i)} ({\Bbb
C})$ and let $\tau \in S(A)$ be a weak limit point of the sequence $\{
\tau_{n(i)} \circ \phi_{i} \} \subset S(A)$.  An easy calculation
shows that $\tau$ is a tracial state.  $\Box$

One should not be tempted to think that the trace constructed above is
faithful.  Of course some very nice unital QD $C^*$-algebras, like the
unitization of the compact operators, can't have a faithful tracial
state.  But we do have the following immediate corollary.

\begin{corollary}
Every simple unital QD $C^*$-algebra has a faithful
trace.
\end{corollary}


\section{Easy Functorial Properties}

The following two facts are immediate from the definition.

\begin{proposition}
A subalgebra of a QD $C^*$-algebra is also QD. 
\end{proposition} 

\begin{proposition}
The unitization of a QD $C^*$-algebra is also QD.
\end{proposition} 

We need some notation before going further. 

\begin{definition}
{\em Let $\{ A_n \}$ be a sequence of $C^*$-algebras.  Then $\Pi_{n
\in {\Bbb N}} A_n = \{ (a_{n}) : \sup_{n} \| a_{n} \| < \infty \}$,
where $(a_{n})$ is an element of the set theoretic product of the
$A_{n}$'s.  We let $\oplus_{n \in {\Bbb N}} A_{n}$ denote the ideal of
$\Pi_{n \in {\Bbb N}} A_{n}$ which consists of elements $(a_{n})$ with
the property that $\lim_{n \to \infty} \| a_n \| = 0$.}
\end{definition}

If $A$ and $B$ are QD and $\pi : A \to B(H)$, $\rho : B \to B(K)$ are
faithful representations whose ranges are quasidiagonal sets of
operators then one easily checks that $A \oplus B$ is QD by
considering the representation $\pi \oplus \rho$.  The following fact
is an easy extension of this argument.

\begin{proposition}
The direct product of QD $C^*$-algebras is QD.  That is, if $\{
A_{n} \}$ is a sequence of
$C^*$-algebras then $\Pi_{n \in {\Bbb N}} A_{n}$ is QD if and only if
each $A_{n}$ is QD.
\end{proposition} 

Recall that if $A$ and $B$ are $C^*$-algebras with faithful
 representations $\pi : A \to B(H)$ and $\rho : B \to B(K)$ then the
 minimal (or `spatial') tensor product is defined to be the
 $C^*$-algebra generated by the image of the algebraic tensor product
 representation $\pi \odot \rho : A \odot B \to B(H) \odot B(K)
 \subset B(H\otimes K).$ The following result, which appeared first in
 [Had2], is left as an easy exercise.  The proof only depends on
 the fact that the tensor product of two finite rank projections is
 again a finite rank projection.

\begin{proposition}
The minimal tensor product of QD $C^*$-algebras is again QD.
\end{proposition} 

If both $A$ and $B$ contain projections and $A \otimes_{min} B$ is QD
then both $A$ and $B$ must be QD as well.  But in general the converse
of Proposition 7.5 is not true (since cones and suspensions are always
QD).

When one of the algebras happens to be nuclear then there is only one
possible tensor product and hence quasidiagonality is always preserved
in this case.  In particular this fact implies that quasidiagonality
is even  invariant under the weaker notion of {\em stable
isomorphism} (cf.\ [BrL1]). (Recall that $A$ and $B$ are stably
isomorphic if $A \otimes {\cal K} \cong B \otimes {\cal K}$, where
${\cal K}$ denotes the compact operators on an infinite dimensional
separable Hilbert space.  Recall also that for separable algebras this
is the same as strong Morita equivalence; cf.\ [BGR].)

It is not known whether or not Proposition 7.5 holds for other tensor
products.  In particular the following question is still open.

\begin{question}
If $A$ and $B$ are QD then is $A \otimes_{max} B$ also QD?
\end{question}


\section{Quotients}

We already pointed out in Example 3.16 that every $C^*$-algebra is a
quotient of a QD $C^*$-algebra.  Thus quasidiagonality does not pass
to quotients in general.  In this section we give a sufficient
condition for a quotient of a QD algebra to be QD.  However, this
condition is far from necessary and it is not clear what the real
obstruction is.  Also, at the end of this section we give an
elementary proof of Corollary 5.3 (i.e.\ one which does not depend on
the fact that cones are always QD).

To state our result we first need a definition. The notion of relative
quasidiagonality was introduced by Salinas in connection with
KK-theory (cf.\ [Sa2]).

\begin{definition}
{\em Let $A$ be a $C^*$-algebra with (closed, 2-sided) ideal $I$.
Then $A$ is said to be {\em quasidiagonal relative to $I$} if $I$ has
an approximate unit consisting of projections which is quasicentral in
$A$. }
\end{definition}

\begin{example}
In general, an algebra can be quasidiagonal relative to an ideal
without itself (or the ideal) being QD.  For example, let $\{ A_i
\}_{i \in {\Bbb N}}$ be a sequence of unital (non-QD) $C^*$-algebras.
Then $A = \Pi_i A_i$ is quasidiagonal relative to the ideal $I =
\oplus_i A_i$.  But the terminology is inspired by a close connection
in the case that the ideal is the compact operators.  Indeed, if $B
\subset B(H)$ is a $C^*$-algebra then it is easy to see that $B$
is a quasidiagonal set of operators if and only if $B + {\cal K}(H)$
is quasidiagonal relative to ${\cal K}(H)$ (cf.\ Proposition 3.6).
\end{example}

\begin{proposition}
Assume $A$ is unital, QD, quasidiagonal relative to an ideal $I$ and
$\pi : A \to A/I$ is locally liftable (cf.\ Definition 2.12). Then $A/I$
is also QD.
\end{proposition}
{\it\noindent Proof.} Let ${\cal F} \subset A/I$ be a finite set and
$\varepsilon > 0$.  In the notation of Corollary 4.5 we let $\phi :
X_{{\cal F}{\cal F}} \to A$ be a unital completely positive splitting. 

Now take a quasicentral approximate unit of projections, say $\{
p_{n} \}$ and consider the (isometric -- though no longer unital)
completely positive splittings $\phi_{n} (x) = (1 - p_{n})
\phi(x) (1 - p_{n})$.  We claim that for sufficiently large
$n$, these maps are $\varepsilon$-multiplicative on ${\cal F}$
and hence from Corollary 4.5  we will have that $A/I$ is QD.  

To see the $\varepsilon$-multiplicativity we first recall that if $a
\in A$ and $\dot{a}$ denotes it's image in $A/I$ then $\| \dot{a} \| =
\lim \| (1 - p_{n})a \|$ since $\{ p_{n} \}$ is an
approximate unit for $I$.  However, since the $p_n$'s are
projections and quasicentral for $A$ we see that $\| \dot{a} \| = \lim
\| (1 - p_{n})a (1 - p_{n})\|$ as well.  Now for $a,b \in
A$ consider the following estimates (see also the proof of Lem.\ 3.1
in [Ar] where these estimates are given in greater generality):
\begin{multline*}    
 \| \phi_{n}(\dot{a}\dot{b}) - \phi_{n}( 
     \dot{a})\phi_{n}(\dot{b}) \| \\[2mm]
   \begin{aligned}
          & =  \|(1 - p_{n})\phi(\dot{a}\dot{b})(1 - p_{n}) - 
                 (1 - p_{n})\phi(\dot{a})(1 - p_{n})
                 \phi(\dot{b})(1 - p_{n}) \| \\[2mm] 
          & \leq \| (1 - p_{n})\bigg(\phi(\dot{a}\dot{b}) - 
                 \phi(\dot{a})\phi(\dot{b}) \bigg)(1 - p_{n})\|\\[2mm] 
          & \phantom{\leq} + 
                 \| (1 - p_{n})\bigg( (1 - p_{n})\phi(\dot{a}) 
                 - \phi(\dot{a})(1 - p_{n}) \bigg) \phi(\dot{b}) 
                 (1 - p_{n}) \|. 
  \end{aligned}
\end{multline*} 
Finally since $\phi$ is a splitting, $\| \dot{x} \| = \lim \| (1 -
p_{n})x (1 - p_{n})\|$ and $\{ p_{n} \}$ is quasicentral we see that
$\phi_{n}$ is $\varepsilon$-multiplicative on ${\cal F}$ for
sufficiently large $n$.  $\Box$

\begin{corollary}
If $A$ is unital, locally reflexive (e.g.\ exact or nuclear), QD and
quasidiagonal relative to an ideal $I$ then $A/I$ is also QD.
\end{corollary}
{\it\noindent Proof.}  Use the previous proposition together with
Theorem 2.13.  $\Box$

\begin{remark} 
The proof of Proposition 8.3 given here is simply a formalization of a
well known argument in the case the ideal is the compact
operators. (cf.\ [D\u{a}d1, Prop.\ 4.5].)
\end{remark}

\begin{example}
Proposition 8.3 is no longer true without the `local liftability'
hypothesis.  Indeed, S. Wassermann gave the first examples of
quasidiagonal sets of operators whose image in the Calkin algebra was
a non-QD $C^*$-algebra (cf.\ [Was1,2]).  Hence, by the remarks in
Example 8.2, Wassermann's examples show that the `local liftability'
hypothesis can't be dropped in Proposition 8.3.
\end{example}

In section 10 we will see that the `right' obstruction to look at for
extensions is probably given in K-theoretic terms (for a large class
of algebras).  However, the following example shows that this is not
the case for the quotient question.  Indeed, it is not at all clear
what type of obstruction one should be looking at in relation to the
quotient question.

\begin{example}
Let $A = {\cal O}_2$ be the Cuntz algebra on two generators (cf.\
[Cu]).  Then we have the short exact sequence $0 \to SA \to CA \to A
\to 0$, where $SA$ and $CA$ denote the suspension and cone,
respectively.  The point we wish to make is that any potential
K-theoretic obstruction would vanish for this example since the six
term exact sequence is trivial.  However, $CA$ is QD while $A$ is not.
\end{example}

Finally we give an elementary proof of the fact that every separable
$C^*$-algebra is a quotient of an RFD algebra (cf.\ Corollary 5.3).  We
will need the following noncommutative generalization of the Tietze
Extension Theorem.

\begin{theorem}(cf.\ [We, Thm.\ 2.3.9])
Let $A$ be a separable $C^*$-algebra and $\pi : A \to B$ be a {\em
surjective} *-homomorphism.  Then $\pi$ extends to a surjective
*-homomorphism $\tilde{\pi} : M(A) \to M(B)$ of multipier algebras.
\end{theorem}

\begin{proposition}
Let $R$ be an RFD algebra.  Then the multiplier algebra, $M(R)$, is
also RFD.
\end{proposition}
{\it\noindent Proof.}  Let $\pi_{n} : R \to A_{n}$ be a sequence of
{\em surjective} *-homomorphisms such that each $A_{n}$ is unital and
QD (e.g.\ finite dimensional) and the map $\oplus_{n \in {\Bbb N}}
\pi_{n}$ is faithful.  Construct extending morphisms $\tilde{\pi}_{n}
: M(R) \to A_{n}$. (In this case the extensions are easy to construct.
For each $n$ let $e_n \in R$ be a lift of the unit of $A_n$ and simply
define $\tilde{\pi}_{n} (x) = \pi_n (xe_n)$ for all $x \in M(R)$.)

In general, if $I \subset A$ is an essential ideal and $\phi : A \to
B$ is a *-homomorphism such that $\phi|_I$ is injective then $\phi$
must be injective on all of $A$ (since any nonzero ideal of $A$ must
have nonzero intersection with $I$).  Hence we see that the
*-homomorphism $\oplus_{n} \tilde{\pi}_{n} : M(R) \to \Pi A_{n}$ is
also injective.  $\Box$

\begin{corollary}
Let $H$ be a (separable) Hilbert space.  Then $B(H)$ is a quotient of a
RFD algebra.
\end{corollary}
{\it\noindent Proof.} Since the compact operators on a separable
Hilbert space are a quotient of an RFD algebra (note that for QD
$C^*$-algebras we do not have to use cones in the proof of Corollary
5.3), it follows from the Tietze Extension Theorem and proposition 8.9
that $B(H) = M({\cal K})$ is a quotient of an RFD algebra. $\Box$


\section{Inductive Limits}

It follows easily from Theorem 3.11  that an
inductive limit of QD $C^*$-algebras where the connecting maps are all
{\em injective} will again be a QD $C^*$-algebra.  We will see
that, in general, inductive limits of QD algebras need not be QD.
However, if the algebras in the sequence are also locally reflexive
(e.g.\ exact or nuclear) then the limit algebra must be QD.

In [BK1] the notion of MF $C^*$-algebra was introduced.  There are a
number of characterizations of these algebras and hence we can choose
the most convenient as our definition (though it is actually a
theorem).

\begin{definition}(cf.\ [BK1, Thm.\ 3.2.2])
{\em A $C^*$-algebra $A$ is MF if and only if $A$ is isomorphic to a
subalgebra of $\Pi M_{n(i)}({\Bbb C}) / \oplus M_{n(i)}({\Bbb C})$ for
some sequence $\{ n_{(i)} \}$.}
\end{definition}

\begin{proposition} (cf.\ [BK1, Prop.\ 3.1.3])
Let $C \subset B(H)$ be a $C^*$-algebra which is also a quasidiagonal
set of operators and $\pi : B(H) \to Q(H)$ denote the quotient map
onto the Calkin algebra.  Then $\pi(C)$ is MF.
\end{proposition}
{\noindent\it Proof.}  By Theorem 5.2 we can find a block diagonal
algebra $B \subset B(H)$ such that $C + {\cal K} = B + {\cal K}$.  If
$P_1 \leq P_2 \leq P_3 \leq \ldots$ are finite rank projetions which
commute with $B$ and converge to $1_H$ then there is a canonical
identification $\Pi_i M_{n(i)}({\Bbb C}) \hookrightarrow B(H) =
B(\oplus_{i \in {\Bbb N}} {\Bbb C}^{n(i)})$, where $n(i) = rank(P_i) -
rank(P_{i - 1})$ (and $P_0 = 0$).  Note that under this identification
we have $\Pi M_{n(i)}({\Bbb C}) \cap {\cal K} = \oplus M_{n(i)}({\Bbb
C})$ and $B \subset \Pi M_{n(i)}({\Bbb C})$.  Hence $$\pi(C) \cong
B/(B \cap {\cal K}) \hookrightarrow \Pi M_{n(i)}({\Bbb C}) / \oplus
M_{n(i)}({\Bbb C}). \ \ \ \ \ \ \ \ \ \Box$$

Evidently every MF algebra has such an extension by the compacts.

It follows that every QD $C^*$-algebra is MF (cf.\ Theorem 3.11).  The
converse is not true by Example 8.6.

The following simple result shows that MF algebras can also be
described as the class of $C^*$-algebras arising as inductive limits
of RFD $C^*$-algebras.

\begin{proposition}
A $C^*$-algebra $A$ is MF if and only if $A$ is isomorphic to an
inductive limit of RFD algebras.
\end{proposition} 
{\it\noindent Proof.}  That an inductive limit of MF algebras (e.g.\
RFD algebras) is again MF is a bit out of the scope of this article.
Please see [BK1, Cor.\ 3.4.4] for the proof.  We will prove the
converse however.

By pulling back the embedding $A \hookrightarrow \Pi M_{n(i)}({\Bbb
C}) / \oplus M_{n(i)}({\Bbb C})$ we can find an RFD algebra $R$ with
ideal $I = \oplus M_{n_{(i)}}({\Bbb C})$ such that $A \cong R/I$.
Consider the finite dimensional ideals $I_k = \oplus_{i = 1}^{k}
M_{n_{(i)}}({\Bbb C})$.  Evidently $R/I_k$ is again an RFD algebra (being
a direct summand of $R$) and hence the natural inductive
system $$R \to R/I_1 \to R/I_2 \to R/I_3 \to \cdots $$ consists of RFD
algebras.  Moreover, since $I = \overline{\cup I_k}$ it is routine to
verify that $A \cong R/I$ is isomorphic to the inductive limit of the
above sequence.  $\Box$

\begin{remark}
It follows that inductive limits of QD $C^*$-algebras {\em need not}
be QD. To get such examples, we let $A \subset B(H)$ be a
$C^*$-algebra which is a quasidiagonal set of operators and such that
the image, $B \subset Q(H)$, in the Calkin algebra is non-QD. (We
mentioned in Example 8.6 that Wassermann  has constructed such
algebras.)  By Propositions 9.2 and 9.3, $B$ is an inductive limit
of RFD algebras which is not QD.
\end{remark}

In contrast to the previous remark, the next result shows
that  mild assumptions will ensure the quasidiagonality of the
limit.  

\begin{theorem}
Let $\{ A_m, \phi_{n,m} \}_{m \in {\Bbb N}}$ be an inductive system of
unital locally reflexive  QD $C^*$-algebras with limit $A =
\lim\limits_{\to} A_i$.  Then $A$ is QD.
\end{theorem}
{\it\noindent Proof.} To clarify our notation, we mean that for all $n
\geq m$ there is a *-homomorphism $\phi_{n,m} : A_m \to A_n$ and we
have the usual compatibility condition that $\phi_{n,m} \circ
\phi_{m,l} = \phi_{n,l}$ whenever $l \leq m \leq n$. We also let
$\Phi_n : A_n \to A$ denote the induced *-homomorphism.

Unitizing the inductive system, if necessary, we may assume that all
the connecting maps are unit preserving.  Now let $\Psi_m : A_m \to
\Pi A_i$ be the *-monomorphism defined by
$$\Psi_m (x) = 0 \oplus \cdots \oplus 0 \oplus x \oplus \phi_{m + 1,
m} (x) \oplus \cdots$$ and $B = C^*( \cup \Psi_m (A_m)) + \oplus A_i
\subset \Pi A_i$.  Then it is easy to see that $B$ is QD,
quasidiagonal relative to the ideal $\oplus A_i$ and $A \cong
B/(\oplus A_i)$. Thus it suffices to see (by Proposition 8.3 and
Remark 4.6) that the quotient map $B \to B/(\oplus A_i)$ is locally
liftable on a dense set.  But this follows from the fact that each
$A_n$ is locally reflexive (cf.\ Theorem 2.13), the maps $\Psi_n$ are
injective, the exact sequences $0 \to (\Psi_n(A_n) \cap \oplus A_i)
\to \Psi_n(A_n) \to \Phi_n(A_n) \to 0$ and the fact that the union of
the $\Phi_n(A_n)$'s is dense in $A$.  $\Box$

\begin{remark} 
Blackadar and Kirchberg have shown that generalized inductive limits
(where the connecting maps are completely positive contractions) of
nuclear QD algebras are again QD (cf.\ [BK1, Cor.\ 5.3.5]).
\end{remark}

Inductive limit decompositions have played a crucial role in (the
finite case of) Elliott's Classification Program.  The next result of
Blackadar and Kirchberg may turn out to have important consequences in
this program.  This theorem follows immediately from [BK1, Prop.\
6.1.6] and [BK2, Cor.\ 5.1].

\begin{theorem}
Let $A$ be a unital simple nuclear QD $C^*$-algebra. Then $A =
\overline{\cup R_i}$ where $R_i \subset R_{i + 1}$ are nuclear RFD
algebras.
\end{theorem}

The remarkable point of this theorem is that the connecting maps in
the inductive system are all {\em injective}.  Indeed, if one relaxes
this condition then we can easily get every nuclear QD $C^*$-algebra.

\begin{proposition}
Let $A$ be a nuclear $C^*$-algebra.  Then $A$ is QD if and only if $A$
is isomorphic to an inductive limit of nuclear RFD $C^*$-algebras.
\end{proposition}
{\noindent\it Proof.} $(\Longleftarrow)$ This follows from Theorem 9.5. 

$(\Longrightarrow)$ This follows from the proof of Proposition 9.3
since extensions of nuclear $C^*$-algebras are nuclear.  $\Box$


\section{Extensions}

Since the Toeplitz algebra is an extension of the compacts by $C({\Bbb
T})$, it follows that extensions of QD algebras need not be QD.
Indeed, as with the quotient question, the general extension problem
for QD algebras appears to be very hard.  As we will see, it is not
even clear whether or not a split extension of QD algebras should be
QD.

We begin, however, with two simple positive results.  The first states
that if the {\em ideal} is sufficiently quasidiagonal then the middle
algebra is always QD.  The second states that if the {\em extension}
is sufficiently quasidiagonal then the middle algebra is always QD.

\begin{proposition}
Assume $0 \to I \stackrel{\iota}{\to} E \stackrel{\pi}{\to} B \to 0$
is exact with $I$ an RFD algebra and $B$ a QD algebra.
Then $E$ is QD.
\end{proposition}
{\it\noindent Proof.}  Let $\phi : E \to M(I)$ be the natural
extension of the inclusion $I \hookrightarrow M(I)$ (cf.\ [We,
2.2.14]).  Then the map $\phi \oplus \pi : E \to M(I) \oplus B$ is
injective.  But Proposition 8.9 states that $M(I)$ is QD (even RFD)
and hence $E$ is a QD $C^*$-algebra.  $\Box$

Note that the proof of Proposition 8.9  actually shows that $M(I)$ is QD
whenever $I$ has a separating family of {\em unital} QD
quotients. Hence the proposition above remains true for ideals of the
form $I = R\otimes_{min} B$ where $R$ is RFD and $B$ is unital and QD.
Hence a natural question is the following.

\begin{question}
Which (nonunital) $C^*$-algebras have QD multiplier algebras?
\end{question}

\begin{definition}
{\em Let $0 \to I \stackrel{\iota}{\to} E \stackrel{\pi}{\to} B \to 0$
be a short exact sequence of $C^*$-algebras.  Such a sequence is
called a {\em quasidiagonal extension} if $E$ is quasidiagonal
relative to $\iota (I)$ (cf.\ Definition 8.1).}
\end{definition}

\begin{remark}
It is important to note that in general an extension being
quasidiagonal has nothing to do with whether or not the middle algebra
$E$ is QD (see Example 8.2).
\end{remark}

\begin{proposition}
Let $0 \to I \stackrel{\iota}{\to} E \stackrel{\pi}{\to} B \to 0$ be a
quasidiagonal extension where both $I$ and $B$ are QD.  Then $E$ is
QD.
\end{proposition}
{\it\noindent Proof.}  To ease notation somewhat, we identify $I$ with
$\iota (I)$ and let $\{ P_{n} \} \subset I$ be an approximate unit of
projections which is quasicentral in $E$.  Now consider the
contractive completely positive maps $\phi_{n} : E \to I \oplus B$,
$\phi_{n} (x) = P_{n}xP_{n} \oplus \pi(x)$.  Evidently these maps are
asymptotically multiplicative.  So we may appeal to Corollary 4.5 and
deduce that $E$ is QD as soon as we verify the following assertion:

{\it Claim.} If $x \in E$ then $\| x \| = max\{ \ \liminf_{n} \|
P_{n}xP_{n} \|, \ \| \pi(x) \| \ \}$.

To prove the claim we pass to the double dual $E^{**}$.  Let $P \in
I^{**} \subset E^{**}$ be the (weak) limit of the $P_{n}$'s.  Then $P$
is central in $E^{**}$ and we have a decomposition $E^{**} = I^{**}
\oplus B^{**}$. Hence (regarding $E \subset E^{**}$) for each $x \in
E$ we have $\| x \| = max\{ \ \| PxP \|, \ \| (1 - P)x(1 - P) \| \
\}$.  But $\| \pi(x) \| = \| (1 - P)x(1 - P) \|$ and $\| PxP \| \leq
\liminf_{n} \| P_{n}xP_{n} \|$ since $P_{n}xP_{n} \to PxP$ in the
strong operator topology. But this proves the claim since the
inequality $\| x \| \geq max\{ \ \liminf_{n} \| P_{n}xP_{n} \|, \ \|
\pi(x) \| \ \}$ is obvious.  $\Box$

\begin{remark}
As mentioned previously, Propositions 10.1 and 10.5 can be regarded as
saying that quasidiagonality is always preserved, provided that either
the {\em ideal} or the {\em extension} is sufficiently quasidiagonal.
This is not true if only the {\em quotient} is highly QD (e.g.\ the
Toeplitz algebra).  Instead a K-theoretic obstruction appears to
govern in general.
\end{remark}

We would now like to discuss the general question of when
quasidiagonality is preserved in extensions.  However, to illustrate
the difficulty of this problem we first pose two basic (open)
questions.

\begin{question}
Let $0 \to I \stackrel{\iota}{\to} E \stackrel{\pi}{\to} B \to 0$ be a
{\em split} exact sequence (i.e.\ there exists a *-homomorphism $\rho
: B \to E$ such that $\pi \circ \rho = id_B$) with $I$ and $B$ QD.  Is
$E$ necessarily QD?
\end{question}

\begin{question}
Let $I$ and $B$ be QD $C^*$-algebras and $\pi : B \to M( I
\otimes {\cal K})$ be a *-monomorphism such that $\pi(B) \cap (I \otimes
{\cal K}) = \{ 0 \}$.  Is $\pi (B) + I \otimes {\cal K}$ necessarily
QD?
\end{question}

Clearly an affirmative answer to Question 10.7 would imply an
affirmative answer to Question 10.8.  In fact the converse is true.

\begin{lemma}
Questions 10.7 and 10.8 are equivalent.
\end{lemma}
{\it\noindent Proof.}  Assume Question 10.8 has an affirmative answer
and let $0 \to I \stackrel{\iota}{\to} E \stackrel{\pi}{\to} B \to 0$
be a split exact sequence and $\rho : B \to E$ be such that $\pi \circ
\rho = id_B$.  Identify $I$ with $\iota (I)$. Let $\eta : E \to B(H)$
be a faithful essential representation.  Then from Theorem 3.11, $\eta
(I) + {\cal K}$ is QD.  Moreover, $\eta (I) + {\cal K}$ is an
essential ideal in $\eta (E) + {\cal K}$.  So replacing $I$ by $\eta
(I) + {\cal K}$ and $E$ by $\eta (E) + {\cal K}$ we may further assume
that $I$ is essential in $E$.  But then $0 \to I \otimes {\cal K}
\stackrel{\iota}{\to} E \otimes {\cal K} \stackrel{\pi}{\to} B \otimes
{\cal K} \to 0$ is still a split exact sequence with $I \otimes {\cal
K}$ essential in $E \otimes {\cal K}$.  Hence $E\otimes {\cal K}$ may
be regarded as a subalgebra of $M(I \otimes {\cal K})$ (cf.\ [We,
2.2.14]) and thus an affirmative answer to Question 10.8 would imply
that $E\otimes {\cal K}$ is QD.  $\Box$

In [BND] it is shown that Question 10.8 has an affirmative answer
under the additional hypothesis that either $I$ or $B$ is
nuclear. Note, however, that even in the case that $I = {\Bbb C}$,
Question 10.8 is not trivial (an affirmative answer still depends on
the full power of Voiculescu's Theorem; cf.\ Theorem 3.11).  Hence it
is not clear whether or not we should expect an affirmative answer to
these questions in general.

If we restrict to the class of nuclear $C^*$-algebras then some
progress can be made on the general extension problem.  Blackadar and
Kirchberg have asked whether or not every nuclear stably finite
$C^*$-algebra is QD (cf.\ [BK1, Question 7.3.1]).  Hence one may ask
whether the extension problem can be solved for stably finite
$C^*$-algebras.  J. Spielberg has given a complete answer to this
question in his work on the AF embeddability of extensions of
$C^*$-algebras.

\begin{proposition} [Sp, Lem.\ 1.5]
Let $0 \to I \to E \to B \to 0$ be an exact sequence with both $I$ and
$B$ stably finite.  If $\partial : K_1 (B) \to K_0 (I)$ denotes the
boundary map of this sequence then $E$ is stably finite if and only if
$\partial(K_1 (B)) \cap K_0^+ (I) = \{ 0 \}$, where $K_0^+ (I)$ is the
canonical positive cone of $K_0 (I)$.
\end{proposition}  

Though the proof is fairly straightforward, we will not prove this
result here as we do not wish to introduce the K-theory which is
needed.

In light of the previous proposition and the question of whether or
not the notions of quasidiagonality and stable finiteness coincide in
the class of nuclear $C^*$-algebras, the following question becomes
quite natural. 

\begin{question}
Let $0 \to I \to E \to B \to 0$ be an exact sequence with both $I$ and
$B$ nuclear QD $C^*$-algebras.  Is it true that $E$ is QD if and only
if $\partial(K_1 (B)) \cap K_0^+ (I) = \{ 0 \}$?
\end{question}

If one approaches this problem via KK-theory then it is probably
necessary to further assume that $B$ satisfies the Universal
Coefficient Theorem (UCT) of Rosenberg and Schochet (cf.\ [RS]).  In
[BND] it is shown that this question is equivalent to some very
natural questions concerning the K-theory of nuclear QD
$C^*$-algebras.  Moreover, it seems likely that an affirmative answer
to the question above could have important consequences in the
classification program (specifically to the classification of Lin's
TAF algebras; [Li1,2]).

In [BND] we also give a partial solution to the question above. The
techniques used to prove the following result are similar to those from
[Sp]. (See also [ELP] for the case that the quotient is AF.)

\begin{theorem} [BND]
Let $0 \to I \to E \to B \to 0$ be an exact sequence with $I$ QD and
$B$ nuclear, QD and satisfying the UCT.  If $\partial : K_1(B) \to K_0
(I)$ is the zero map then $E$ is QD
\end{theorem}


\section{Crossed Products}

In this section we discuss when crossed products of QD $C^*$-algebras
are again QD.  This is not always the case since the (purely infinite)
Cuntz algebras are stably isomorphic to crossed products of AF
algebras by $\Bbb Z$.  The basic theory of crossed products by locally
compact groups can be found in [Pe1, Chpt.\ 7].  (See also [Dav,
Chpt.\ 8] for a nice treatment of the discrete case.)

We begin with a corollary of an imprimitivity theorem of P. Green.  To
state the result we will need to introduce some notation.  So, let $G$
be a separable locally compact group and $H \subset G$ be a closed
subgroup.  Then $G/H$ (the space of left cosets) is a separable
locally compact space.  There is a natural action $\gamma$ of $G$ on
$C_0 (G/H)$ defined by $\gamma_g (f)(xH) = f(g^{-1}xH)$ for all $xH \in
G/H$ and $f \in C_0 (G/H)$.  The crossed products below are the full
crossed products and all groups actions $\alpha : G \to Aut(A)$ are
assumed to be suitably continuous (i.e.\ for each $a \in A$ the map $g
\mapsto \alpha_g (a)$ is continuous).

\begin{theorem}([Gr2, Cor.\ 2.8 ]) Let $\alpha : G \to Aut(A)$ be a
homomorphism from the separable locally compact group $G$.  For each
closed subgroup $H \subset G$ there is an isomorphism $$A\otimes C_0
(G/H) \rtimes_{\alpha \otimes \gamma} G \cong (A \rtimes_{\alpha|_{H}}
H) \otimes {\cal K},$$ where ${\cal K}$ denotes the compact operators
on a separable (finite dimensional if and only if $G/H$ is finite)
Hilbert space.
\end{theorem}

For the rest of this section we will only be dealing with amenable
groups (cf.\ [Pe1, 7.3]) and hence we do not need to distinguish
between reduced and full crossed products (cf.\ [Pe1, Thm.\ 7.7.7]).

\begin{corollary}
Let $A$ be QD and $\alpha : G \to Aut(A)$ be a homomorphism with $G$ a
separable compact group.  Then $A \rtimes_{\alpha} G$ is QD.
\end{corollary}
{\it\noindent Proof.}  Let $H \subset G$ be the zero subgroup.  The
previous theorem then asserts that $A\otimes C (G) \rtimes_{\alpha
\otimes \gamma} G \cong A \otimes {\cal K}$.  But $A \otimes {\cal K}$
is QD  and there is a natural embedding $A \rtimes_{\alpha} G
\hookrightarrow A\otimes C (G) \rtimes_{\alpha \otimes \gamma} G$
since $G$ amenable implies that the full and reduced crossed products
are isomorphic (cf.\ [Pe1, 7.7.7 and 7.7.9]).  $\Box$

For non-compact discrete groups the problem is considerably harder.
However, Rosenberg has shown that we must restrict to the class of
amenable groups.

\begin{theorem} [Ros, Thm.\ A1]
If $G$ is discrete and $C^*_r (G)$ is QD then $G$ is amenable.
\end{theorem}
 
It is not known whether the converse of this theorem holds (cf.\ [Vo4,
3.1]), but F$\o$lner's characterization of amenable groups in terms of
almost shift invariant finite subsets leads one to believe that the
converse should be true.

For actions of ${\Bbb Z}$ there are only two classes of $C^*$-algebras
where we currently have complete information on the quasidiagonality
of $A \rtimes_{\alpha} {\Bbb Z}$; when $A$ is abelian or AF.  Before
stating the theorems we first give a definition. 

\begin{definition}
{\em Let $A$ be a $C^*$-algebra.  Then $A$ is called {\em AF
embeddable} if there exists a *-monomorphism $\rho : A \to B$ where
$B$ is AF.}
\end{definition}

Of course AF embeddable $C^*$-algebras are QD.  However, it is a
nontrivial fact that the converse is not true. In fact, even RFD
algebras need not be AF embeddable.  The best known example is the
full group $C^*$-algebra $C^* ({\Bbb F}_2)$.  This is RFD but is not
exact and hence cannot be embed into {\em any} nuclear (in particular,
AF) algebra (cf.\ [Was3]). However, for crossed products of abelian or
AF algebras by ${\Bbb Z}$, quasidiagonality does imply AF
embeddability.
\begin{theorem}([Pi, Thm.\ 9])
Let $\phi : X \to X$ be a homeomorphism of the compact metric space
$X$ and $\Phi \in Aut( C(X))$ denote the induced automorphism.  Then
the following are equivalent:
\begin{enumerate}
  \item $C(X) \rtimes_{\Phi} {\Bbb Z}$ is AF embeddable, 
  \item $C(X) \rtimes_{\Phi} {\Bbb Z}$ is QD, 
  \item $C(X) \rtimes_{\Phi} {\Bbb Z}$ is stably finite, 
  \item `$\phi$ compresses no open sets.' (That is,
  if $U \subset X$ is open and $\phi ( \overline{U} ) \subset U$ then
  $\phi (U) = U$.)
\end{enumerate}
\end{theorem}

\begin{theorem}([BrN1, Thm.\ 0.2])
Let $A$ be AF and $\alpha \in Aut(A)$ be given.  Then the following
are equivalent:
\begin{enumerate}
  \item $A \rtimes_{\alpha} {\Bbb Z}$ is AF embeddable, 
  \item $A \rtimes_{\alpha} {\Bbb Z}$ is QD, 
  \item $A \rtimes_{\alpha} {\Bbb Z}$ is stably finite, 
  \item `$\alpha_{*} : K_0 (A) \to K_0 (A)$ compresses no elements.' 
  (That is, if $x \in K_0 (A)$ and $\alpha_{*} (x) \leq x$ in 
   the natural order then $\alpha_{*} (x) = x$.)
\end{enumerate}
\end{theorem}

We have chosen to formulate the above results in a way that
illustrates their similarities.  In both cases the hard implications
are $4 \Rightarrow 1$.  Also in both cases it is not at all clear that
the techniques in the proof will be of much use in general.  Before
going beyond actions of ${\Bbb Z}$ we wish to point out that there is
no harm in assuming unital algebras.

\begin{proposition} 
Let $A$ be nonunital, $\alpha \in Aut(A)$, $\tilde{A}$ be the
unitization of $A$ and $\tilde{\alpha} \in Aut(\tilde{A})$ the unique
unital extension of $\alpha$.  Then $A \rtimes_{\alpha} {\Bbb Z}$ is
QD if and only if $\tilde{A} \rtimes_{\tilde{\alpha}} {\Bbb Z}$ is QD.
\end{proposition}
{\it\noindent Proof.} Recall that we always have a {\em split} exact
sequence $$0 \to A \rtimes_{\alpha} {\Bbb Z} \to \tilde{A}
\rtimes_{\tilde{\alpha}} {\Bbb Z} \to C({\Bbb T}) \to 0.$$ Thus the
implication $(\Leftarrow$) is immediate and ($\Rightarrow$) follows
from Theorem 10.12 since abelian algebras are nuclear, QD and satisfy the
Universal Coefficient Theorem (cf.\ [RS]).  $\Box$

Another natural direction to consider would be to try crossed products
of well behaved $C^*$-algebras by more general groups.  (We must stay
within the class of amenable groups, though, because of Rosenberg's
result; cf.\ Theorem 11.3)  However, even for actions of ${\Bbb Z}^2$ this
is a problem.  Indeed the following question of Voiculescu remains
open even now -- more than 15 years after Pimsner's result for $C(X)
\rtimes_{\Phi} {\Bbb Z}$.

\begin{question}(cf.\ [Vo4, 4.6]) 
When is $C(X) \rtimes_{\Phi} {\Bbb Z}^2$ AF embeddable?
\end{question}

For crossed products of certain simple AF algebras the question is
more manageable. 

\begin{theorem} [BrN2, Thm.\ 1]
If $A$ is UHF and $\alpha : {\Bbb Z}^n \to Aut(A)$ is a homomorphism
then there always exists a *-monomorphism $\rho : A\rtimes_{\alpha}
{\Bbb Z}^n \to B$ where $B$ is AF.
\end{theorem} 

The proof of this result (and Theorem 11.6 above) depends in an
essential way on a technical notion known as the {\em Rohlin property}
for automorphisms.  This notion has been used by Connes, Kishimoto,
Evans, Nakamura and others (with great success!) in classifying
automorphisms of operator algebras.  Moreover, Kishimoto has used
these ideas to prove that many crossed products of certain simple
$A{\Bbb T}$ algebras by automorphisms with the Rohlin property will
again be $A{\Bbb T}$ (which is much stronger than just saying they are
QD).  See, for example, [Kis1-4].

\begin{remark}
One nice consequence of Green's theorem (Theorem 11.1) is that
understanding crossed products by ${\Bbb Z}^n$ gives results about
much more general groups.  For example, if $G$ is a finitely generated
discrete abelian group then $G \cong {\Bbb Z}^n \oplus F$ where $F$ is
a finite (hence compact) abelian group then by Green's result we have
an embedding $A \rtimes_{\alpha} G \hookrightarrow (A
\rtimes_{\alpha|_{{\Bbb Z}^n}} {\Bbb Z}^n)\otimes {\cal K}$.  Writing
a general discrete abelian group as an inductive limit of finitely
generated such groups one can then handle crossed products by
arbitrary discrete abelian groups.  One can then proceed to take
extensions by arbitrary separable compact groups and build a very
large class of groups for which it suffices to consider crossed
products by ${\Bbb Z}^n$.  (See Def.\ 3.4 and the proof of Thm.\ 2 in
[BrN2] for more details).
\end{remark}


\section{Relationship with Nuclearity}  

It was an open question for quite some time whether or not
quasidiagonality implied nuclearity.  In [Had2], Hadwin asked whether
or not every `strongly' quasidiagonal (e.g.\ simple QD) $C^*$-algebra
was nuclear.  Then in [Po], Popa asked whether every simple unital QD
$C^*$-algebra with `sufficiently many projections' (e.g.\ real rank
zero) was nuclear.  There was some evidence supporting a positive
answer to these questions.  The strongest was the following theorem of
Popa.

\begin{theorem} [Po, Thm.\ 1.2]
Let $A$ be a simple unital $C^*$-algebra with `sufficiently many
projections' (e.g.\ real rank zero).  Then $A$ is QD if and only if
for each finite set ${\cal F} \subset A$ and $\varepsilon > 0$ there
exists a (non-zero) finite dimensional subalgebra $B \subset A$ with
unit $P = 1_B$ such that $\| [a, P] \| \leq \varepsilon$ for all $a
\in {\cal F}$ and $P {\cal F} P \subset^{\epsilon} B$ (cf.\ Definition
4.4).
\end{theorem}

The necessity of the technical condition above is quite hard, however
the sufficiency is easily seen.  Indeed, if one assumes the technical
condition then we can find a sequence of finite dimensional
subalgebras $B_n \subset A$ with units $P_n$ such that $\| [a, P_n] \|
\to 0$ and $d(P_n aP_n , B) \to 0$ for all $a \in A$.  Now let $\Phi_n
: A \to B_n$ be a conditional expectation and consider the maps
$\phi_n : A \to B_n$ defined by $\phi_n (a) = \Phi_n (P_n a P_n)$.
This sequence of maps is evidently asymptotically multiplicative and
hence defines a *-homomorphism $$A \to \Pi B_n / \oplus B_n.$$ Since
$A$ is unital this morphism is nonzero and since $A$ is simple, this
morphism is injective.  Hence the maps $\phi_n$ are also
asymptotically isometric which implies (by Theorem 4.2) that $A$ is QD.
(Note that the hypothesis of a unit can't be dropped here.  Indeed,
the stabilization ${\cal K} \otimes A$ of any unital $C^*$-algebra $A$
satisfies the technical condition stated above. Simply take $B_n$ of
the form $C_n \otimes 1_A$ where $C_n$ is almost orthogonal to a large
part of ${\cal K}$.)  

The above result gave one hope of deducing nuclearity via the
characterization in terms of injective enveloping von Neumann algebras
(cf.\ [CE1]). However, it turns out that this is not possible as the
following result of D\u{a}d\u{a}rlat shows.   

\begin{theorem} [D\u{a}d2, Prop.\ 9]
There exists a unital, separable, simple, QD $C^*$-algebra with real
rank zero, stable rank one and unique tracial state which is not exact
(and hence not nuclear).
\end{theorem}

The converse of the question we have been considering above is also
interesting and worth discussion.  Namely, what sort of general
conditions on a $C^*$-algebra imply quasidiagonality?

\begin{example}
A Cuntz algebra ${\cal O}_n$ (cf.\ [Cu]) is simple, separable,
unital, nuclear, has real rank zero and is not QD (since it is purely
infinite; cf.\ Proposition 3.19).
\end{example}

To get a finite non-QD example is a bit more delicate.  Recall that
$C^*_r (G)$, where $G$ is a discrete group, is always stably finite
since it has a faithful tracial state.  Also recall that Rosenberg has
shown that if the reduced group $C^*$-algebra of a discrete group is
QD then the group must be amenable (cf.\ Theorem 11.3).

\begin{example}
Let ${\Bbb F}_2$ denote the free group on two generators.  Then
$C^{*}_{r} ({\Bbb F}_2)$ is simple, unital, separable, exact, has
stable rank one (cf.\ [DHR]) and a unique tracial state but is not QD
since ${\Bbb F}_2$ is not amenable.
\end{example}
 
To get an example with the added property of real rank zero one can
simply consider $C^{*}_{r} ({\Bbb F}_2) \otimes U$, where $U$ is some
UHF algebra (cf.\ [R$\o$r, Thm.\ 7.2]).

It is also interesting to note that there are no known examples of
finite {\em nuclear} non-QD $C^*$-algebras. (Recall that $C^{*}_{r}
({\Bbb F}_2)$ is only exact.)  In fact, as noted in Section 10,
Blackadar and Kirchberg have formulated the following question.

\begin{question} [BK1, Question 7.3.1] 
If $A$ is nuclear and stably finite then must $A$ necessarily be QD?
\end{question}

This question is of particular interest in Elliott's classification
program (cf.\ [Ell]).  Indeed, if this question turns out to have an
affirmative answer then classifying simple unital nuclear {\em finite}
$C^*$-algebras may be equivalent to classifying simple unital nuclear
QD $C^*$-algebras.  (One would still have to resolve the important
open question of whether every simple finite algebra is stably finite
- which is equivalent to the open question of whether every simple
infinite $C^*$-algebra is purely infinite.)  The point is that for
simple QD algebras (with enough projections) one has the structure
theorem of Popa to work with.  In fact, Lin has introduced a class of
$C^*$-algebras (the so-called TAF algebras; [Li1]) whose definition
looks similar to Popa's structure theorem.  Moreover, there are
classification results for some of these TAF algebras (cf.\ [Li2],
[DE1]) and it is not unreasonable to think that someday the general QD
case can be handled in ways similar to the current strategies being
applied to the TAF case.


\section{More Advanced Topics}

In our final section we will present some miscellaneous results which
don't quite fit into any of the previous sections.  The first is a
very important result of Voiculescu which shows that quasidiagonality
is a homotopy invariant.  Recall that two $C^*$-algebras $A$ and
$B$ are called {\em homotopic} if there exist *-homomorphisms $\phi :
A \to B$ and $\psi : B \to A$ such that $\phi \circ \psi$ is homotopic
to $id_B$ and $\psi \circ \phi$ is homotopic to $id_A$ (cf.\
[Bl2], [We]).  

\begin{theorem}
Let $A$ and $B$ be homotopic $C^*$-algebras.  Then $A$ is QD if and
only if $B$ is QD.
\end{theorem}

Voiculescu actually proved a more general result (cf.\ [Vo3, Thm.\
5]).  In [D\u{a}d1, Thm.\ 1.1] D\u{a}d\u{a}rlat generalized this to show that
quasidiagonality is even an invariant of the weaker notion of
`asymptotic completely positive homotopy equivalence'.  As 
mentioned previously,  this result implies that the cone over {\em
any} $C^*$-algebra is QD since cones are homotopic to $\{ 0 \}$.

Free products of $C^*$-algebras were introduced in [Av] and
independently in [Vo5]. (See also [VDN].)  Reduced free products are
rarely QD.  The standard example of a reduced free product is $C^*_r
({\Bbb F}_2) = C^* ({\Bbb Z}) \ast C^* ({\Bbb Z})$, where the reduced
free product is taken with respect to Haar measure on the circle.  The next
result of F. Boca is in stark contrast. (See also [ExLo] where the
class of RFD algebras is shown to be closed under full free products.)

\begin{theorem} [Bo, Prop.\ 13]
If $A$ and $B$ are unital QD $C^*$-algebras, then the full free
product (amalgamating over the units) $A \ast B$ is also QD.
\end{theorem}

We next point out the connection between quasidiagonality and the
notions of projectivity and semiprojectivity.  These notions are
studied at length in [Lo].

\begin{definition}
{\em Let $A$ be a $C^*$-algebra.  Then $A$ is called {\em projective}
if for every $C^*$-algebra $B$, closed 2-sided ideal $I \subset B$ and
*-homomorphism $\phi : A \to B/I$ there exists a lifting
*-homomorphism $\psi : A \to B$.  $A$ is called {\em semiprojective}
if for every $C^*$-algebra $B$, closed 2-sided ideal $I \subset B$
such that $I = \overline{\cup_n I_n}$ for ideals $I_1 \subset I_2
\subset \ldots$ and *-homomorphism $\phi : A \to B/I$ there exists an
$n$ and a lifting *-homomorphism $\psi : A \to B/I_n$ (that is, a
lifting for the canonical quotient map $B/I_n \to B/I$).  }
\end{definition}

The projective case in our next result is well known.  The
semiprojective case was pointed out by B. Blackadar, though his proof
was different.

\begin{proposition}
If $A$ is projective then $A$ is RFD.  If $A$ is MF and semiprojective
then $A$ is RFD.
\end{proposition}
{\noindent\it Proof.}  First assume that $A$ is projective.  By
Corollary 5.3 $A$ is a quotient of an RFD algebra.  But then the
definition of projectivity implies that $A$ embeds into an RFD algebra
and hence is itself RFD.

Now assume that $A$ is semiprojective and MF.  By the proof of
Proposition 9.3 we can find an RFD algebra $R$ with finite dimensional
ideals $I_n \subset I_{n + 1}$ such that $A \cong R/I$ where $I =
\overline{\cup_n I_n}$.  The definition of semiprojectivity then
provides an embedding $A \hookrightarrow R/I_n \subset R$ for some
$n$.  $\Box$

We now discuss a beautiful connection between quasidiagonality and the
question of whether or not `Ext is a group'. (See also the discussion
in [Vo4].)  Here we mean the classical BDF Ext semigroups.  Recall
that if $A$ is nuclear then the Choi-Effros lifting theorem implies
that $Ext(A)$ is a group. (See [Ar] for a very nice treatment of this
theory.)  But it is known that there exist $C^*$-algebras $A$ for
which $Ext(A)$ is not a group (cf.\ [An], [Was1,2]).  However these
examples are not ``natural'' examples (though those in [Was1] are
quotients of ``natural'' examples).  Indeed, it has been a long
standing open problem to determine whether or not $Ext(C^{*}_{r}
({\Bbb F}_2))$ is a group.  It is believed that $Ext(C^{*}_{r}
({\Bbb F}_2))$ is not a group and we now outline one approach to
proving this.

We described the class of MF algebras in Section 9.  Recall that
these algebras can be characterized as those which appear as the image
in the Calkin algebra of a quasidiagonal set of operators in $B(H)$
(cf.\ Proposition 9.2).

\begin{corollary}
Let $A$ be MF and assume $Ext(A)$ is a group.  Then $A$ is QD.
\end{corollary}
{\it\noindent Proof.}  If $Ext(A)$ is a group then every
*-monomorphism $\phi : A \to B(H)/{\cal K}$ has a completely positive
lifting (cf.\ [Ar, pg.\ 353]).  But then from Propositions 8.3  we
see that $A$ must be QD.  $\Box$

It follows then that every nuclear MF algebra is QD.  Recall, though,
that there exist non-QD MF algebras. But it is not known whether
Wassermann's examples are exact. The following question remains open.

\begin{question}
Do there exist exact non-QD MF algebras?  In particular is $C^{*}_{r}
({\Bbb F}_2)$ MF?
\end{question}

Kirchberg has also proved some remarkable results connecting
quasidiagonality, $Ext$ and various lifting properties of
$C^*$-algebras (see [Kir1]).

Finally, we wish to point out a connection with one of the most
important questions in $C^*$-algebras.  Namely, whether or not the
Universal Coefficient Theorem (UCT) holds for all nuclear separable
$C^*$-algebras (cf.\ [RS]).  We will not formulate this question
precisely as it is well out of the scope of these notes.  However, the
experts will have no problem following our argument.  The main
ingredient is the following `two out of three principle' for the UCT.

\begin{theorem} (cf.\ [RS, Prop.\ 2.3 and Thm.\ 4.1])
Let $0 \to I \to E \to B \to 0$ be a short exact sequence with $E$
nuclear and separable.  If any two of $\{ I, \ E, \ B \}$ satisfy the
UCT then so does the third.  In particular, if $I$ and $E$ satisfy the
UCT then so does $B$.
\end{theorem}

Our final result has been noticed by several experts.

\begin{corollary}
If the UCT holds for all separable nuclear RFD algebras then the UCT
holds for all separable nuclear $C^*$-algebras.
\end{corollary}
{\it\noindent Proof.}  By the two out of three principle, it suffices
to show that every separable nuclear $C^*$-algebra is a quotient of a
separable nuclear RFD algebra.  But this is contained in
Corollary 5.3  $\Box$


\section{Further Reading}

Below are references to some of the topics around quasidiagonality
which are only briefly discussed (or not discussed at all) in these
notes.

{\noindent\it AF embeddability.} [BrN1,2], [D\u{a}d4,5], [Li3], [Pi],
[PV2], [Sp], [Vo2].

{\noindent\it Ext and KK-theory.} [BrL2], [DE2], [DHS], [Kir1], [PV1],
[Sa1,2], [Wa1,2].

{\noindent\it Classification.}  [DE1], [Ell], [Li2] and their bibliographies.

{\noindent\it MF, (strong) NF algebras and inner quasidiagonality.} [BK1,2].

{\noindent\it General.} [Had2], [Th], [Vo4].


\section{Appendix: Nonseparable QD $C^*$-algebras}

In this appendix we treat the case of nonseparable $C^*$-algebras.
Hence we no longer require the Hilbert spaces in this section to be
separable either.  The results of this section (in particular
Corollary 15.7) are necessary for the general case of Voiculescu's
characterization of QD $C^*$-algebras.  Though we have seen some of
these results stated in the literature, we have been unable to find
any proofs and hence complete proofs will be given.

\begin{definition}
{\em A subset $\Omega \subset B(H)$ is a called a {\em quasidiagonal
set of operators} if for each finite set $\omega \subset \Omega$,
finite set $\chi \subset H$ and $\varepsilon > 0$ there exists a
finite rank projection $P \in B(H)$ such that $\| [T, P] \| \leq
\varepsilon$ and $\| P(x) - x \| \leq \varepsilon$ for all $T \in
\omega$ and $x \in \chi$.}
\end{definition}

It is still easy to see that a set $\Omega \subset B(H)$ is a
quasidiagonal set of operators if and only if the $C^*$-algebra
generated by $\Omega$, $C^*(\Omega) \subset B(H)$, is a quasidiagonal
set of operators.

We may finally give the general definition of a quasidiagonal
$C^*$-algebra.

\begin{definition}
{\em Let $A$ be a $C^*$-algebra.  Then $A$ is called {\em
quasidiagonal} (QD) if there exists a faithful representation $\pi : A
\to B(H)$ such that $\pi(A)$ is a quasidiagonal set of operators.}
\end{definition}

There is one subtle point that needs resolved here.  Namely we must
show that the previous definition is equivalent to Definition 3.8 in
the case that $A$ is a separable $C^*$-algebra. 

\begin{lemma}
Let $A$ be a separable $C^*$-algebra and assume that there exists a
faithful representation $\pi : A \to B(H)$ such that $\pi(A)$ is a
quasidiagonal set of operators.  Then there exists a faithful
representation $\rho : A \to B(K)$ such that $K$ is a separable
Hilbert space and $\rho(A)$ is a quasidiagonal set of operators.
\end{lemma}
{\it\noindent Proof.} Let $\pi : A \to B(H)$ be a faithful
representation such that $\pi(A)$ is a quasidiagonal set of
operators. We will show that there exists a {\em separable} subspace
$K \subset H$ which is $\pi(A)$-invariant and such that the
restriction representation $\rho = \pi_K = P_K \pi(\cdot) P_K : A \to
B(K)$ (cf.\ Definition 3.9) is faithful and has the property that
$\rho(A)$ is a quasidiagonal set of operators.

The idea is to construct an increasing sequence of separable
$\pi(A)$-invariant subspaces $\tilde{K}_1 \subset \tilde{K}_2 \subset
\tilde{K}_3 \ldots$ and finite rank projections $Q_n$ such that $Q_n
(H) \subset \tilde{K}_{n + 1}$, $\| [Q_n, \pi(a) ] \| \to 0$ for all
$a \in A$ and $\| Q_n (\xi) - \xi \| \to 0$ for all $\xi \in \cup
\tilde{K}_i$.  If we further arrange that the restriction of $\pi(A)$
to $\tilde{K}_1$ is faithful then it is clear that $K = \overline{\cup
\tilde{K}_i}$ is the desired subspace.
  
We begin by choosing a sequence $\{ a_i \} \subset A$ which is dense
in the unit ball of $A$.  For each $n \in {\Bbb N}$ we then choose a
sequence of unit vectors $\{ \xi^{(n)}_i \}_{i \in {\Bbb N}} \subset
H$ such that $\| \pi(a_i)\xi^{(n)}_i \| > \| a_i \| - 1/2^n$. Let $K_1
\subset H$ be the closure of the span of $\{ \xi^{(n)}_i \}_{i,n \in
{\Bbb N}}$ and let $\tilde{K}_1$ be the closure of $\pi(A)K_1$.  Then
it is clear that $\tilde{K}_1$ is separable, $\pi(A)$-invariant and
the restriction of $\pi(A)$ to $\tilde{K}_1$ is faithful (since it is
isometric on $\{ a_i \}$).

One then constructs the desired $\tilde{K}_i$ and $Q_i$ recursively as
follows.  Let $\{ h^{(1)}_{i} \}$ be an orthonormal basis for
$\tilde{K}_1$.  Choose a finite rank projection $Q_1 \in B(H)$ such
that $\| [Q_1, \pi(a_1) ] \| < 1/2$ and $Q_1 (h^{(1)}_{1}) =
h^{(1)}_{1}$. Recall from the proof of Proposition 3.4 that we can
always arrange the stronger condition $Q_1 (h^{(1)}_{1}) =
h^{(1)}_{1}$.

Next let $X_2 = span\{ Q_1(H), \tilde{K}_1 \}$, $\tilde{K}_2$ be the
closure of $\pi(A) X_2$ and let $\{ h^{(2)}_{i} \}$ be an orthonormal
basis for $\tilde{K}_2$.  Now choose a finite rank projection $Q_2 \in
B(H)$ such that $\| [Q_2, \pi(a_i) ] \| < 1/(2^2)$ for $i = 1,2$, $Q_2
(h^{(j)}_{i}) = h^{(j)}_{i}$ for $i,j = 1,2$ and $Q_1 \leq Q_2$ (this
is arranged by requiring that $Q_2 (h) = h$ for a (finite) basis of
$Q_1 (H)$).

Next let $X_3 = span\{ Q_2(H), \tilde{K}_2 \}$, $\tilde{K}_3$ be the
closure of $\pi(A) X_3$ and let $\{ h^{(3)}_{i} \}$ be an orthonormal
basis for $\tilde{K}_3$, etc.  Proceeding in this way we get an
increasing sequence of separable $\pi(A)$-invariant subspaces
$\tilde{K}_1 \subset \tilde{K}_2 \subset \tilde{K}_3 \ldots$ with
orthonormal bases $\{ h^{(n)}_{i} \}_{i \in {\Bbb N}}$ and finite rank
projections $Q_n \leq Q_{n + 1}$ such that $Q_n (h^{(j)}_{i}) =
h^{(j)}_{i} $ for $i,j = 1, \ldots, n$ and $Q_n (H) \subset
\tilde{K}_{n + 1}$, $\| [Q_n, \pi(a) ] \| \to 0$ for all $a \in A$.
Evidently this proves the lemma.  $\Box$

Hence we see that Definitions 3.8 and 14.2 are equivalent for separable
$C^*$-algebras.  Indeed, it clear that if $A$ is separable and
satisfies Definition 3.8 then $A$ also satisfies Definition 14.2.  On the
other hand, if $A$ is separable and satisfies Definition 14.2 then by
the previous lemma we can find a representation of $A$ on a
{\em separable} Hilbert space which gives a quasidiagonal set of
operators and hence $A$ satisfies Definition 3.8 as well.

We will need the following elementary, but technical, lemma.

\begin{lemma}
Let $\pi : A \to B(H)$ be a faithful representation where $A$ is
separable (but $H$ is not).  Then there exists a separable
$\pi(A)$-invariant subspace $K \subset H$ with the property that
$\pi_K : A \to B(K)$ is faithful, $\pi_K (a)$ is a finite rank
operator if and only if $\pi(a)$ is a finite rank operator and in this
case $dim(\pi(a)H) = dim(\pi_K (a)K)$.
\end{lemma}
{\noindent\it Proof.}  The idea is to find a sequence of
$\pi(A)$-invariant separable subspaces, $\tilde{H}_i$, with the
following properties:
\begin{enumerate}
   \item The restriction of $\pi(A)$ to $\tilde{H}_1$ is faithful.
   \item If $a \in A$ is such that $\pi(a)$ is a finite rank operator 
         then $\pi(a)H \subset  \tilde{H}_1$.
   \item $\tilde{H}_m \ \mbox{\Large $\perp$} \ 
         (\tilde{H}_1 \oplus \tilde{H}_2 \oplus
         \ldots \oplus \tilde{H}_{m - 1})$ 
   \item If $P_{\tilde{H}_m} \pi(a) P_{\tilde{H}_m} = 0$ then $(1 -
         P_{\tilde{H}_1 \oplus \ldots \oplus \tilde{H}_{m - 1}})
         \pi(a) (1 - P_{\tilde{H}_1 \oplus \ldots \oplus \tilde{H}_{m
         - 1}}) = 0$ for all $a \in A$, where for any subspace 
         $L \subset H$, $P_L$ denotes the orthogonal projection onto $L$.
\end{enumerate} 
Having the subspaces $\{ \tilde{H}_i \}$ we define $$K = \oplus_{i =
1}^{\infty} \tilde{H}_i \subset H$$ and note that $\pi_K : A \to B(K)$
is faithful (since this was already arranged on $\tilde{H}_1$).
Moreover, condition 2 ensures that if $\pi(a)$ is a finite rank
operator then $dim(\pi(a)H) = dim(\pi_K (a)K)$.  Finally, note that if
$\pi_K (a)$ is a finite rank operator then there exists some integer
$m \in {\Bbb N}$ such that $\tilde{\pi}_m (a) = P_{\tilde{H}_m} \pi(a)
P_{\tilde{H}_m} = 0$.  Hence
\begin{multline*}
  \begin{aligned} 
     \pi(a) 
       &  =  \tilde{\pi}_1 (a) \! \oplus \! \cdots \! \oplus \! 
           \tilde{\pi}_{m - 1} (a) \! \oplus \! (1 \! - \! 
           P_{\tilde{H}_1 \oplus \ldots
           \oplus \tilde{H}_{m - 1}}) \pi(a) (1 \! - \! P_{\tilde{H}_1 \oplus
           \ldots \oplus \tilde{H}_{m - 1}}) \\[2mm]
       & = \tilde{\pi}_1 (a) \oplus \cdots \oplus
           \tilde{\pi}_{m - 1} (a) \oplus 0 \\[2mm] 
       & = \pi_K (a),
  \end{aligned}
\end{multline*}
by condition 4 above.  Hence $\pi(a)$ is also a finite rank operator
and clearly $dim(\pi(a)H) = dim(\pi_K (a)K)$.  So we now show how to
construct subspaces $\tilde{H}_i$ as above.

Begin by letting ${\cal F}(A) = \{ a \in A : dim(\pi(a)H) < \infty \}$
and choosing a countable dense subset $\{ a_i \}_{i \in {\Bbb N}}
\subset {\cal F}(A)$.  For each $i \in {\Bbb N}$ let $L_i = \pi(a_i)H$
and define $H_1$ to be the closure of $$span\{ \bigcup_{i =
1}^{\infty} \pi(A)L_i \}.$$ By throwing in a countable number of
vectors (as in the proof of Lemma 15.3) we can replace $H_1$ with a
larger $\pi(A)$-invariant subspace $\tilde{H}_1$ such that the
restriction of $\pi(A)$ to $\tilde{H}_1$ is also faithful.  We claim
that this $\tilde{H}_1$ also satisfies condition 2 above.  Indeed, if
$a \in {\cal F}(A)$ then we can find a subsequence $a_{i_j} \to a$.
But since $\pi(a_{i_j})H \subset \tilde{H}_1$ and $\pi(a_{i_j}) \to
\pi(a)$ it is clear that $\pi(a)H \subset \tilde{H}_1$ as well.  Hence
we have constructed $\tilde{H}_1$ with the desired properties.

Assume now that we have constructed $\tilde{H}_1, \ldots ,
\tilde{H}_{m - 1}$ with the desired properties.  To get $\tilde{H}_m$
we simply consider the separable $C^*$-algebra $$C = (1 -
P_{\tilde{H}_1 \oplus \ldots \oplus \tilde{H}_{m - 1}}) \pi(A) (1 -
P_{\tilde{H}_1 \oplus \ldots \oplus \tilde{H}_{m - 1}}).$$ By the
proof of Lemma 15.3 we can find a separable $C$-invariant subspace
$\tilde{H}_m \subset (1 - P_{\tilde{H}_1 \oplus \ldots \oplus
\tilde{H}_{m - 1}})H$ such that the restriction of $C$ to
$\tilde{H}_m$ is faithful.  Evidently $\tilde{H}_m$ is also $\pi(A)$
invariant, perpendicular to $\tilde{H}_1 \oplus \tilde{H}_2 \oplus
\ldots \oplus \tilde{H}_{m - 1}$ and condition 4 above is nothing more
than the statement that the map $C \to P_{\tilde{H}_m}C
P_{\tilde{H}_m}$ is faithful.  $\Box$

As in section 3 we want to resolve the technical issue of
nondegeneracy of representations.

\begin{lemma}
Let $A$ be a $C^*$-algebra and $\pi: A \to B(H)$ be a faithful
representation.  Let $L \subset H$ be the nondegeneracy subspace of
$\pi(A)$ and $\pi_L : A \to B(L)$ denote the restriction.  Then
$\pi(A)$ is a quasidiagonal set of operators if and only if $\pi_L
(A)$ is a quasidiagonal set of operators.
\end{lemma}
{\it\noindent Proof.}  The implication $( \Leftarrow )$ is proved
exactly as in Lemma 3.10.  Also, if $A$ is unital, the implication
$(\Rightarrow)$ is the same and so we only have to show
$(\Rightarrow)$ in the case that $A$ is nonunital.

So assume that $A$ is nonunital and $\pi(A)$ is a quasidiagonal set of
operators.  Note that we cannot apply Voiculescu's Theorem in this
setting since the dimensions of $H$ and $L$ may be different.  To
resolve this problem we first note that since quasidiagonality is
defined via finite sets it suffices to show that $\pi_L (B)$ is a
quasidiagonal set of operators for every {\em separable}
$C^*$-subalgebra $B \subset A$.  

Given a finite set of vectors $\chi \subset L$, by Lemma 15.4, we can
find a separable subspace $K \subset L$ with the property that $\chi
\subset K$, $K$ is $\pi_L (A)$-invariant, the restriction to $K$ is
faithful, $\pi_L (a)$ is finite rank if and only if $\pi_K (a)$ is
finite rank and in this case $rank( \pi_L (a)) = rank( \pi_K (a))$.
As in the proof of Lemma 15.3 we can now enlarge $K$ to a separable
$\pi(A)$-invariant subspace $\tilde{K} \subset H$ (we do not have
$\tilde{K} \subset L$, of course) such that $\pi_{\tilde{K}} (A)$ is a
quasidiagonal set of operators.  Since we have been careful about
separability and preservation of rank it now follows from Voiculescu's
theorem (version 2.6) that $\pi_{\tilde{K}}$ and $\pi_K$ are
approximately unitarily equivalent and hence $\pi_K(A)$ is a
quasidiagonal set of operators.  $\Box$

\begin{theorem}
Let $\pi : A \to B(H)$ be a faithful {\em essential} (cf.\ Definition
2.8) representation.  Then $A$ is QD if and only if $\pi (A)$ is a
quasidiagonal set of operators.
\end{theorem}
{\it\noindent Proof.} Clearly we only have to prove the necessity. As
in the proof of the previous lemma, it suffices to show that $\pi(B)$
is a quasidiagonal set of operators for every separable subalgebra $B
\subset A$.  

Let $\chi \subset H$ be an arbitrary finite set and use Lemma 15.4 to
construct a separable $\pi(B)$-invariant subspace $K \subset H$ such
that $\chi \subset K$ and the restriction to $K$ is both faithful and
essential.  The remainder of the proof is now similar to that of
Theorem 3.11.  $\Box$

The next corollary shows that with  care, one can usually just
treat the separable case when dealing with quasidiagonality.

\begin{corollary}
$A$ is QD if and only if all of it's finitely generated subalgebras
are QD.
\end{corollary}
{\it\noindent Proof.}  The necessity is obvious from the definition.
So assume all finitely generated subalgebras of $A$ are QD and let
$\pi : A \to B(H)$ be a faithful essential representation.  Then for
each finitely generated subalgebra $B \subset A$ the restriction
$\pi|_B$ is a faithful essential representation and hence (by Theorem
15.6) $\pi (B)$ is a quasidiagonal set of operators.  It then follows
that $\pi(A)$ is a quasidiagonal set of operators. $\Box$

Finally we observe the nonseparable version of Theorem 4.2. 

\begin{corollary}[Voiculescu]
Let $A$ be a $C^*$-algebra.  Then $A$ is QD if and only if $A$
satisfies {\em ($\ast$)}.
\end{corollary}
{\it\noindent Proof.}  As in the proof of Theorem 4.2, we may assume
that $A$ is unital.  From Arveson's Extension Theorem it follows that
$A$ satisfies ($\ast$) if and only if every separable unital
subalgebra of $A$ satisfies ($\ast$).  Similarly, from Corollary 15.7
it follows that $A$ is QD if and only if every separable unital
subalgebra of $A$ is QD.  Hence this corollary follows from the
separable case.  $\Box$

\vspace{3mm}

{\noindent\bf Acknowledgements.}  These notes are based on a series of
lectures given by the author at the University of Tokyo in the spring
of 1999.  We gratefully acknowledge the support of the NSF, which
allowed us to spend one year in Tokyo, and express our thanks to the
seminar participants for enduring our lectures.  Special thanks also
go out to my thesis advisor, Marius D\u{a}d\u{a}rlat, and Larry Brown for
teaching me so much about QD $C^*$-algebras.


\end{document}